\pdfoutput=1
\def\REVIEW{0}

\ifdefined\REVIEW
    \documentclass[
        a4paper,authoryear,12pt
    ]{article}
    \usepackage{lineno, setspace}

\else
    \documentclass[
        a4paper,
        twocolumn,authoryear,12pt
    ]{article}
\fi

\usepackage{IEK10} 
\usepackage{natbib}

\usepackage{placeins}

\def\IEK10{
  Forschungszentrum Jülich GmbH,
  Institute of Energy and Climate Research,
  Energy Systems Engineering (IEK-10),
  Jülich 52425,
  Germany
}
\def\RWTH{
  RWTH Aachen University
  Aachen 52062,
  Germany
}
\def\JARA{
  JARA-ENERGY,
  Aachen 52056,
  Germany
}
\def\SVT{
  RWTH Aachen University,
  Process Systems Engineering (AVT.SVT),
  Aachen 52074,
  Germany
}
\def\LTT{
  RWTH Aachen University,
  Institute of Technical Thermodynamics (LTT),
  Aachen 52062,
  Germany
}

\newcommand{\mytitle}{Optimal design of a local renewable electricity supply system for power-intensive production processes with demand response}

\newcommand{\affil}{
  \begin{itemize}[leftmargin=3mm, itemsep=0mm]
    \item[$^a$]\IEK10
    \item[$^b$]\RWTH
    \item[$^c$]\LTT
    \item[$^d$]\JARA
    \item[$^e$]\SVT
  \end{itemize}
}

\newcommand{\myauthor}{Sonja H. M. Germscheid$^{a,b}$\orcidlink{0000-0002-0411-7711}, 
Benedikt Nilges$^{c}$\orcidlink{0000-0002-9510-2039},
Niklas von der Assen$^{c,d}$\orcidlink{0000-0001-8855-9420},
Alexander Mitsos$^{d,a,e}$\orcidlink{0000-0003-0335-6566}, Manuel Dahmen$^{a,*}$\orcidlink{0000-0003-2757-5253}}
\newcommand{\correspondingauthor}{M. Dahmen, \IEK10 \newline E-mail: m.dahmen@fz-juelich.de}
\author{\myauthor}

\usepackage[  hidelinks]{hyperref}
\usepackage[capitalise, nameinlink]{cleveref}
\crefname{table}{Tab.}{Tab.}

\usepackage{orcidlink}  
\usepackage{cancel}
\usepackage[export]{adjustbox} 
\usepackage{graphicx}
\usepackage{subcaption}

\usepackage[acronym]{glossaries}

\usepackage{pgfplots}
\pgfplotsset{compat=newest}
\usepgfplotslibrary{groupplots}
\usepgfplotslibrary{dateplot}
\begin{document}

\newacronym{da}{DA}{day-ahead}
\newcommand{\da}{\gls*{da}}
\newacronym{id}{ID}{intraday}
\newcommand{\id}{\gls*{id}}
\newacronym{dr}{DR}{demand response}
\newcommand{\dr}{\gls*{dr}}
\newacronym{tac}{TAC}{total annualized cost}
\newcommand{\tac}{\gls*{tac}}
\newacronym{gwi}{GWI}{global warming impact}
\newcommand{\gwi}{\gls*{gwi}}

\newcommand{\maxP}{\ensuremath{\theta_\text{max}}}
\newcommand{\Pnom}{\ensuremath{P_\text{nom}}}
\newcommand{\dt}{\ensuremath{\Delta{t}}}
\newcommand{\nS}{\ensuremath{\mathbb{S}}}

\ifx\REVIEW\undefined
\twocolumn[
\begin{@twocolumnfalse}
\fi
  \thispagestyle{firststyle}
  \begin{center}
    \begin{large}
      \textbf{\mytitle}
    \end{large} \\
    \myauthor
  \end{center}
  \vspace{0.5cm}
  \begin{footnotesize}
    \affil
  \end{footnotesize}
  \vspace{0.5cm}
  \doublespacing
  \begin{abstract}
    This work studies synergies arising from combining industrial demand response and local renewable electricity supply. 
    To this end, we optimize the design of a local electricity generation and storage system with an integrated demand response scheduling of a continuous  power-intensive production process in a multi-stage problem. 
    We optimize both total annualized cost and global warming impact and consider local photovoltaic and wind electricity generation, an electric battery, and electricity trading on day-ahead and intraday market. 
    We find that installing a battery can reduce emissions and enable large trading volumes on the electricity markets, but significantly increases cost.
    Economically and ecologically-optimal operation of the process and battery are driven primarily by the electricity price and grid emission factor, respectively,  rather than locally generated electricity. 
    A parameter study reveals that cost savings from the local system and flexibilizing the process behave almost additively.
    \end{abstract}

\vspace{0.5cm}

\noindent \textbf{Keywords}: \textit{Integrated design and scheduling, stochastic programming, demand response, local electricity supply system, renewable energy
}

  \vspace*{5mm}
\ifx\REVIEW\undefined
\end{@twocolumnfalse}
]
\fi

 \doublespacing
 \newpage
 \section{Introduction}

Renewable electricity has a varying supply that leads to time-varying electricity prices on the electricity markets. The time-varying prices can incentivize flexible industrial processes to adapt their momentary production rate and, thus, power consumption in a \dr{} scheduling, which can reduce operational cost and is considered electricity grid balancing \citep{Daryanian1989,Zhang2016_DSMbasics,Burre2020,Mitsos2018}.
\dr{} savings can be  improved by participating in multiple short-term electricity markets, see, e.g., \cite{Leo2020,DalleAve2019,Simkoff2020,Liu2016,PandZic2013,Kwon2017,Golmohamadi2018,nolzen2022market,
Germscheid_2022,Germscheid2023,Schafer2019_bidding,Varelmann2022}.
Furthermore, flexible operation should be accounted for at design stage in order to determine optimal investment decisions for both the production processes itself, see, e.g., \cite{Mitra2014_ASU_integrated,Teichgraeber2020,Steimel2015,Seo2023,Leenders2019}, and for its local energy supply system, see, e.g., \cite{Yunt2008,Zhang2019_integrated1,Voll2013,Baumgartner2019,Langiu2022_ORC,Bahl2017, Fleschutz2023}.

In local energy supply systems, integrated design and scheduling has already been used to optimize on-site renewable electricity generation and storage systems. In the corresponding studies, the considered systems satisfy a fixed demand profile but can offer flexibility by combining different electricity generation technologies, see, e.g., \cite{Zhang2019_integrated1,Bahl2017, Fleschutz2023,Baumgartner2019}.
Furthermore, combining on-site renewable electricity generation and storage systems with flexible production processes can reduce both production cost and CO$_2$ emissions, which has been shown for, e.g., 
a water electrolyzer in combination with Power-to-X processes \citep{Mucci2023}, ammonia and nitric acid production \citep{Wang2020}, ammonia generation \citep{Allman2018}, and methanol production \citep{Martin2016}.
Further power-intensive, flexible production processes could benefit from the combination with on-site renewable electricity supply, e.g., the chlor-alkali electrolysis \citep{Bree2019}, seawater desalination \citep{Ghobeity2010}, or air separation \citep{Ierapetritou2002}.
However, to the best of our knowledge, a generalized assessment about synergies arising from combining on-site electricity supply systems and \dr{}-capable processes has not been conducted yet.

In our prior work \citep{Germscheid_2022}, we conducted a  \dr{}  potential assessment of power-intensive production processes by means of the generalized process model introduced by \cite{Schafer2020_generic}.
The generalized process model can represent a wide range of continuous production processes by means of few key process characteristics, i.e., oversizing, minimal part load, ramping limitations, and production storage capacity. 
We analyzed the benefit of participating simultaneously in both the \da{} and the \id{} spot electricity market, but neglected potential electricity provision by on-site renewable electricity generation and storage \citep{Germscheid_2022}.

In this article, we extend our prior work \citep{Germscheid_2022} by integrating the scheduling of the generalized production process into the design optimization of a local electricity generation and storage system.
In the resulting multi-stage approach, we optimize the design of the local renewable electricity supply system considering photovoltaic (PV) power,  wind power, and an electric battery for a location in Germany. On the lower stages, we optimize the \dr{}  scheduling of both the energy system and the process, together with the electricity market participation. 
We consider both economic and ecologic design objectives, i.e., we optimize with respect to the \tac{} and the \gwi{}, respectively.
We study the influence of different degrees of process flexibility on the optimal design of the energy system and the resulting ecologic and economic savings.
Similar to our prior work \citep{Germscheid_2022}, we consider simultaneous market participation in both the \da{} and \id{} electricity market  to analyze  the benefit of considering multiple electricity markets in an integrated design and scheduling problem.

The remainder of the article is structured as follows:
\cref{sec:structure} explains the structure of the integrated design and scheduling problem. We specify the objectives in \cref{sec:objective} and the operational constraints in \cref{sec:constraints}. 
The scenarios and the model parameters are specified in \cref{sec:scenarios} and \cref{sec:specifications}, respectively.
We discuss the optimal energy system design for a reference process in \cref{sec:basicDesign}, the dependency between process parameters and potential savings in \cref{sec:synergy}, and the benefit of considering simultaneous \da{} and \id{} market participation at design stage in \cref{sec:simultanousmarket}. In \cref{sec:conclusion}, we conclude our work.

\section{Methods}\label{sec:methods}

\subsection{Structure of the integrated design and scheduling problem} \label{sec:structure}

\begin{figure}[!ht]
\centering
  \includegraphics[trim=0cm 23.4cm 7.5cm 0.5cm ,clip, width =0.9\textwidth]{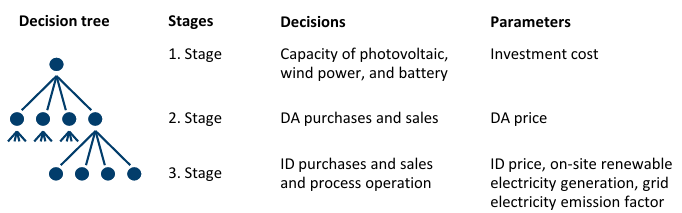}
\caption{Structure of the integrated design and scheduling problem: The problem structure allows optimizing the design of the local electricity generation and storage system while considering simultaneous \da{} and \id{} market participation and process \dr{}. The decision tree models the  chronological sequence of decisions. Each branch represents the realization of an uncertain parameter and each node represents a decision point.
\label{fig:structure}
} 
\end{figure}

Integrated design and scheduling problems are often set up as two-stage stochastic problems \citep{Birge2011} with the design decisions on the first stage and scheduling decisions and operational constraints on the second stage, see, e.g., \cite{Yunt2008,Zhang2019_integrated1,Langiu2022_ORC,Mitra2014_ASU_integrated,Teichgraeber2020,Steimel2015,Seo2023,Bahl2018}.
In this work, we determine the optimal design of a local electricity supply system for a flexible industrial production process.
We account for simultaneous \da{} and \id{} market participation in the design and scheduling problem by the three-stage structure shown in \cref{fig:structure}. In the first stage, the design decisions for the energy system are made, i.e., photovoltaic (PV), wind power, and electric battery capacities are to be determined. 
The \da{} trading decisions are taken the day before the operation when the \id{} price and renewable electricity generation are still uncertain.
Thus, we consider the \da{} decisions on the second stage and \id{} trading and operational decisions on the third stage. In particular, the operation of the flexible process is adapted on the third stage in response to realizations of the \id{} price, renewable electricity production, and the momentary emission factor of the grid electricity. 
Note that we consider a time-varying average grid emission factor similar to \cite{Baumgartner2019} and \cite{Nilges2023} and that the emission factor is uncertain a day before the actual consumption.

Similar to our prior work \citep{Germscheid_2022}, we consider hourly \da{} purchases and quarter-hourly \id{} purchases and sales and assume a one-day scheduling horizon.
In addition, we allow selling electricity from the local generation and storage system on the \da{} market.
For simplicity, we assume that throughout any quarter-hour time slice, renewable electricity generation is constant.

In the following, we omit a distinct notation for second- and third-stage parameters and variables for better readability.
In particular, we consider \da{} trading decisions $\boldsymbol{q}_\text{DA,$s_3$}$ and \da{} price $\boldsymbol{c}_\text{DA,$s_3$}$ on the third stage instead of the second stage and guarantee equality on the second stage by means of non-anticipativity constraints:
\begin{align}
    \boldsymbol{c}_\text{DA,$s_3$} = \boldsymbol{c}_{\text{DA},\hat{s}_3} \quad &\text{if} \quad f(s_3)=f(\hat{s}_3) \quad \forall (s,\hat{s})\in \mathbb{S}_3\times\mathbb{S}_3,\\
    \boldsymbol{q}_\text{DA,$s_3$} = \boldsymbol{q}_{\text{DA},\hat{s}_3} \quad &\text{if} \quad f(s_3)=f(\hat{s}_3) \quad \forall (s,\hat{s})\in \mathbb{S}_3\times\mathbb{S}_3
\end{align}
Here, $f:\mathbb{S}_3\xrightarrow{}\mathbb{S}_2$ maps a node on the third stage of the decision tree, i.e., $s_3\in\mathbb{S}_3$, to the respective node on the second stage, i.e., $s_2\in\mathbb{S}_2$. 
In the following, we refer to $s\in \mathbb{S}$ instead of $s_3\in\mathbb{S}_3$ for conciseness and we use $s$ for indexing the scenarios, i.e., paths in the decision tree.

\subsection{Objectives}\label{sec:objective}

We consider the \acrfull{tac} and the \acrfull{gwi} as economic and ecologic objective, respectively.

The \tac{} is defined as
\begin{align}
&\text{\tac{}} = \text{CAPEX} +\text{OPEX}_{el} +\text{OPEX}_{Grid}, \label{eqn:tac}\\
\text{with} \quad\quad&\text{CAPEX} = \sum_{i\in\{\text{PV,W,B}\}} \left( \frac{(\gamma_\text{1}+1)^{\gamma_{2,i}}\gamma_{1}}{(\gamma_{1}+1)^{\gamma_{2,i}}-1}  \text{CAPEX}^0_i \;Q_i +\gamma_{3,i}\;Q_i \right),  \label{eqn:capex}\\
&\text{OPEX}_{el} =  365 \sum_{s\in \nS } \pi_{s} (4\;\dt \; {\boldsymbol{c}_\text{DA,s}} \cdot \boldsymbol{q}_\text{DA,s}   + \dt \;\boldsymbol{c}_{\text{ID},s} \cdot \boldsymbol{q}_{\text{ID},s}),  \label{eqn:opex_el}\\
   & \text{OPEX}_{\text{Grid}} =  365 \sum_{s\in \nS } \pi_{s}  \sum_{t=1}^{T} \text{OPEX}_{\text{Grid}, s,t}, \label{eqn:ffee1}\\
    & \text{OPEX}_{\text{Grid}, s,t} \geq c_\text{Fee} \; \dt\;  (q_{\text{DA,s,}\lfloor \frac{t-1}{4} \rfloor+1} +q_{\text{ID,}s,t} ), \label{eqn:ffee2}\\
    & \text{OPEX}_{\text{Grid}, s,t} \geq0 . \label{eqn:ffee3}
\end{align}
In \cref{eqn:tac}, we compute the \tac{} as the sum of investment cost, $\text{CAPEX}$, and operational cost, $\text{OPEX}$.
According to current legislation in Germany (Status 2023), the grid fee has to be paid in addition to the market price for electricity removed from the electricity grid \citep{Gridfee}. Thus, we consider both annual operational cost from electricity procurement, $\text{OPEX}_{el}$, as well as annual grid fee cost, $\text{OPEX}_{Grid}$. 
\cref{eqn:capex} specifies the investment cost of the local electricity generation and storage system considering photovoltaic (PV), wind power (W), and electric battery (B). Similar to \cite{Baumgartner2019}, we calculate the annualized CAPEX based on the total investment cost $\text{CAPEX}^0_i$, the present value factor with interest rate $\gamma_{1}$ and life time $\gamma_{2,i}$ \citep{Broverman2010}, a maintenance factor $\gamma_{3,i}$, and the installed capacity of the respective technology $Q_i$.
Note that in contrast to \cite{Baumgartner2019}, we consider a component-specific life time $\gamma_{2,i}$. 
\cref{eqn:opex_el} defines the electricity cost  $\text{OPEX}_{el}$ by the purchases and sales on hourly \da{} and quarter-hourly \id{} electricity market $\boldsymbol{q}_\text{DA,s}$ and $\boldsymbol{q}_{\text{ID},s}$, the \da{} and \id{} electricity price $\boldsymbol{c}_\text{DA,s}$ and $\boldsymbol{c}_{\text{ID},s}$, the time step size $\Delta t=0.25$h, and the probability $\pi_{s} $ of scenario $s$.
In \cref{eqn:ffee1}, the grid cost $\text{OPEX}_{\text{Grid}}$ is derived from the sum of the grid cost $\text{OPEX}_{\text{Grid}, s,t}$ of each scenario $s$ and time step $t$ with a total of 96 time steps for the one-day scheduling horizon, i.e, $T=96$.
\cref{eqn:ffee2} and \cref{eqn:ffee3} constitute lower bounds for $\text{OPEX}_{\text{Grid}, s,t}$, which ensure that $\text{OPEX}_{\text{Grid}, s,t}$ is equal to zero in case of electricity injection into the grid and greater or equal to the grid fee with the grid fee cost $c_\text{Fee}$ in case of electricity removal from the grid.
$\text{OPEX}_{\text{Grid}, s,t}$ is equal to the respective lower limit, i.e., zero or the grid fee, when minimizing the \tac{}.
Note that in $\text{OPEX}_{\text{Grid}, s,t}$,  $q_{\text{ID,}s,t}$ varies quarter-hourly and $q_{\text{DA,}s,t}$ varies hourly.

The expected annual \gwi{} is computed as:
\begin{align}
    \text{\gwi{}} =& 365 \sum_{s\in \nS } \pi_{s}  \sum_{t=1}^{T} \Bigl(  \text{GWI}^\text{el}_{s,t} \; \dt\;  ( q_{\text{DA,s,}\lfloor \frac{t-1}{4} \rfloor+1}   + q_{\text{ID},s,t}) \Bigl)+
    \sum_{i\in\{\text{PV,W,B}\}}\frac{\text{GWI}^\text{i} Q_{\text{i}}}{\gamma_{2,\text{i}}}
\label{eqn:GWI}
\end{align}
\cref{eqn:GWI} considers the quarter-hourly average \gwi{} of the electricity from the grid $\text{GWI}^\text{el}_{s,t}$, the \gwi{} of the installed PV, wind power, and battery capacity, i.e., $\text{GWI}^\text{PV}$, $\text{GWI}^\text{W}$, and $\text{GWI}^\text{B}$, and their respective life time $\gamma_{2,i}$ in years.
Consequently, the total annual \gwi{} depends on hourly and quarter-hourly purchases and sales from the \da{} and \id{} market $q_{\text{DA,s,t}}$ and $q_{\text{ID,s,t}}$, respectively, and the installed PV, wind, and battery capacities, i.e., $\text{Q}^\text{PV}$, $\text{Q}^\text{W}$, and $\text{Q}^\text{B}$, respectively.
Note that we allow for a \gwi{} credit, i.e., negative \gwi{}, in case of electricity sales accounting for an avoided emission burden \citep{LCA_Credit}.

\subsection{Operational constraints}\label{sec:constraints}

In the following, we  shortly describe the generalized process model from our prior work \citep{Schafer2020_generic,Germscheid_2022} and discuss the operational constraints specific to the local energy system and the electricity trading.

The generalized process model \citep{Schafer2020_generic} relies on few key parameters to describe the \dr{} capabilities. In this work, we consider the key characteristics oversizing, minimal part load, product storage capacity with cyclic storage constraints, and ramping limitation. 
Note that without efficiency losses, the production rate of the process scales directly with the process power intake.
For a detailed explanation of the generalized process model, including the model equations, we refer to \cite{Germscheid_2022}.

In the electricity generation and storage system, we consider the \da{} and \id{} purchases and sales
$q_{\text{DA,s,t}}$ and $q_{\text{ID,s,t}}$, respectively, with positive values corresponding to purchases. In addition, we consider that the electricity purchases, PV power $q_{\text{PV,}s,t}$, wind power $q_{\text{W,}s,t}$, and battery charge and discharge $q_{\text{in,}s,t}$ and $q_{\text{out,}s,t}$ are equal to the power intake $p_{s,t}$ of the production process by  the energy balance:
\begin{align}
    p_{s,t} &= q_{\text{DA,s,}\lfloor \frac{t-1}{4} \rfloor+1} +q_{\text{ID,}s,t} + q_{\text{PV,}s,t} + q_{\text{W,}s,t}  -q_{\text{in,}s,t}+q_{\text{out,}s,t} \label{eqn:ebalance}
\end{align}
Additionally, we consider operational constraints similar to \cite{Baumgartner2019}:
\begin{align}
q_{\text{PV,}s,t}&= \bar{q}_{\text{PV,}s,t}\; Q_\text{PV}, \label{eqn:operation:pv}\\
q_{\text{W,}s,t}&= \bar{q}_{\text{W,}s,t} \;Q_\text{W},\label{eqn:operation:wind}\\
0&\leq q_{\text{in,}s,t}\leq Q_\text{B} / \tau,\label{eqn:operation:bat1}\\
0&\leq q_{\text{out,}s,t}\leq Q_\text{B} / \tau,\label{eqn:operation:bat2}\\
0&\leq SOC_{s,t} \leq Q_\text{B},\label{eqn:operation:bat2b}\\
SOC_{s,t+1}&=SOC_{s,t}+(\eta^\text{in}q_{\text{in,}s,t}- \frac{q_{\text{out,}s,t}}{\eta^\text{out}} )\dt,\label{eqn:operation:bat3}\\
SOC_{s,1}&=SOC_{s,T+1}\label{eqn:operation:bat4}
\end{align}
Here, \cref{eqn:operation:pv} and \cref{eqn:operation:wind} define the PV and wind power generation, $q_{\text{PV,}s,t}$ and $q_{\text{W,}s,t}$, by multiplying the relative power output, $\bar{q}_{\text{PV,}s,t}$ and $\bar{q}_{\text{W,}s,t}$, with the installed PV and wind capacity, $Q_\text{PV}$ and $Q_\text{W}$, respectively.
\cref{eqn:operation:bat1} and \cref{eqn:operation:bat2} constrain the battery charge and discharge, i.e., $q_{\text{in,}s,t}$ and $q_{\text{out,}s,t}$, respectively, by the installed battery capacity $Q_\text{B}$ and the allowed charging and discharging rate $\tau$. 
\cref{eqn:operation:bat2b} constrains the state of charge $SOC_{s,t}$ by the installed battery capacity $Q_\text{B}$.
\cref{eqn:operation:bat3} relates charging and discharging with respective efficiency losses $\eta^\text{in}$ and $\eta^\text{out}$ and the state of charge.
Additionally, we consider the cyclic constraint,  \cref{eqn:operation:bat4}, requiring that the state of charge is the same at the beginning and at the end of the scheduling horizon.

We consider trading electricity on both \da{} and \id{} market while making use of the local energy system:
\begin{align}
- q_{\text{PV,}s,t}- q_{\text{W,}s,t}-Q_\text{B} / \tau &\leq q_{\text{DA,s,}\lfloor \frac{t-1}{4} \rfloor+1} ,\label{eqn:trade1}\\
- q_{\text{PV,}s,t}- q_{\text{W,}s,t}-SOC_{s,t} /\dt &\leq q_{\text{DA,s,}\lfloor \frac{t-1}{4} \rfloor+1} , \label{eqn:trade2} \\
- q_{\text{DA,s,}\lfloor \frac{t-1}{4} \rfloor+1}- q_{\text{PV,}s,t}- q_{\text{W,}s,t} -Q_\text{B} / \tau &\leq q_{\text{ID},s,t} , \label{eqn:trade3}\\
- q_{\text{DA,s,}\lfloor \frac{t-1}{4} \rfloor+1}- q_{\text{PV,}s,t}- q_{\text{W,}s,t} -SOC_{s,t} / \dt &\leq q_{\text{ID},s,t},  \label{eqn:trade4}\\
q_{\text{DA,s,}\lfloor \frac{t-1}{4} \rfloor+1} &\leq \Pnom(1+\maxP)  +Q_\text{B} / \tau ,   \label{eqn:trade5}\\
q_{\text{DA,s,}\lfloor \frac{t-1}{4} \rfloor+1} &\leq \Pnom(1+\maxP)  +\frac{Q_\text{B}-SOC_{s,t}}{\dt}  ,  \label{eqn:trade6}\\
q_{\text{ID},s,t} &\leq \Pnom(1+\maxP) +Q_\text{B} / \tau  , \label{eqn:trade7}\\
q_{\text{ID},s,t} &\leq \Pnom(1+\maxP) + \frac{Q_\text{B}-SOC_{s,t}}{\dt}  \label{eqn:trade8}
\end{align}
Here, \cref{eqn:trade1,eqn:trade2,eqn:trade3,eqn:trade4} constrain \da{} and \id{} sales by produced PV and wind power, the current state of charge, and maximum discharging capabilities of the installed battery.
Similar to our prior work \citep{Germscheid_2022}, \cref{eqn:trade3,eqn:trade4} allow selling previously purchased \da{} electricity on the \id{} market.
\cref{eqn:trade5,eqn:trade6,eqn:trade7,eqn:trade8} constrain \da{} and \id{} purchases by the maximum power consumption of the production process and the battery. The maximum consumption of the process is defined by the nominal consumption \Pnom{} and the process oversizing \maxP{}. The maximum consumption of the battery is given by the state of charge and the charging capabilities of the battery.

\subsection{Time series}\label{sec:scenarios}

In the following, we specify the time series representing realizations of uncertain parameters in the assessment and refer to them as scenarios in this context.

\begin{figure}[!ht]
\begin{subfigure}[t]{.45\textwidth}
    \centering
  \includegraphics[width =\linewidth]{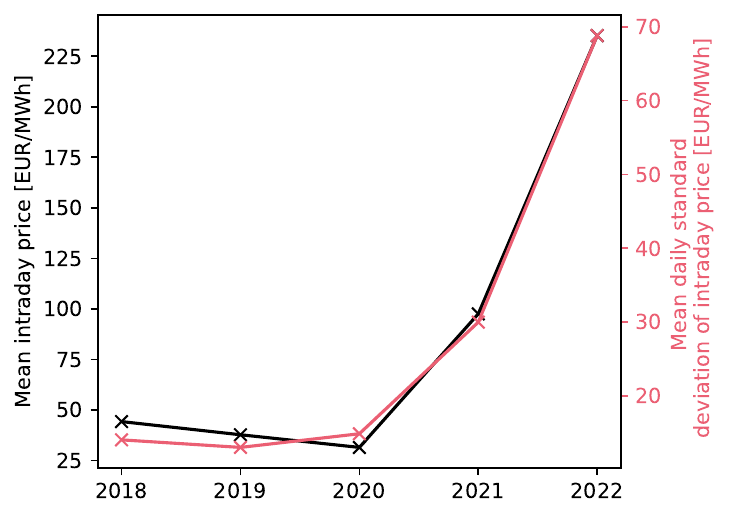}
  \caption{\id{} electricity price }
  \label{fig:ID_price}
  \end{subfigure}
  \hfill
\begin{subfigure}[t]{.45\textwidth}
    \centering
  \includegraphics[width =\linewidth]{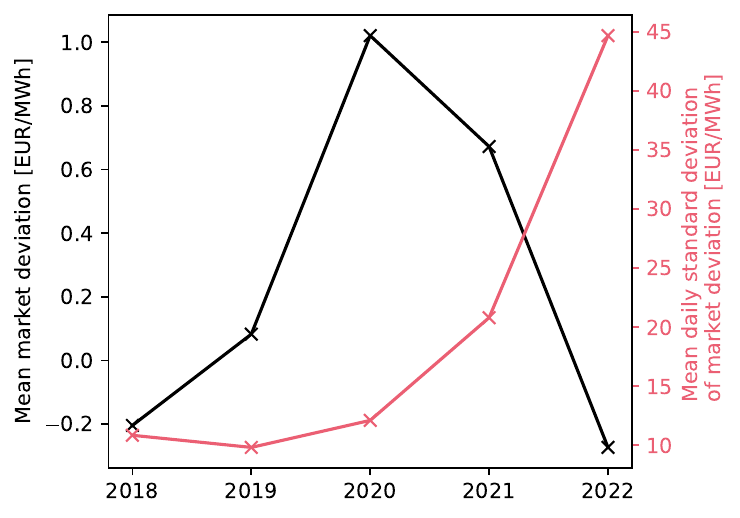}
  \caption{Market deviation }
  \label{fig:Dev_price}
  \end{subfigure}
  \caption{Mean and mean daily standard deviation of \id{} electricity price (a) and market deviation (b), i.e., the price difference between the \da{} and the \id{} price.}
  \label{fig:price_develop}
\end{figure}

We base our scenarios for the electricity price, wind power generation, PV power
generation, and grid emission factor on historical time series. Specifically, 
for the electricity price, we use data from the German \da{} and \id{} spot market \citep{FrauenhoferInstitute2019}, assuming that the consumer can purchase and sell electricity at \da{} market-clearing price and \id{} index price. We refer to \cite{Germscheid_2022} for a detailed explanation about these assumptions.
\cref{fig:price_develop} shows the annual mean and the annual mean daily standard deviation of the \id{} electricity market price and the market deviation, i.e., the difference between the \da{} and the \id{} price.
The corresponding figure for the \da{} price shows similar characteristics as the one of the \id{} price and can be found in Section 1 of the supporting material.
\cref{fig:ID_price} reveals a price decrease in 2020 that can be attributed to the initial phase of the COVID-19 pandemic \citep{Halbrugge2021_covid} and an increase of mean and standard deviation in 2021 and 2022 due to the conflict in the Ukraine \citep{Haucap2022_ukraine}.
\cref{fig:Dev_price} reveals that in 2020 and 2021, the mean \id{} price was slightly larger than the mean \da{} price as indicated by the positive market deviation. Moreover, the standard deviation of the market deviation has significantly increased in 2021 and 2022.

\dr{}  scheduling optimization necessitates  electricity price time series rather than an average annual electricity price. Long-term German electricity price forecasting is challenging, e.g., due to the conflict in Ukraine and the German energy transition. 
The time series for 2030 derived by the project MONA 2030 \citep{MONA} used in \cite{Schafer2020_generic} and the prices for 2050 reported in \cite{DENA} are outdated.  
For our assessment, we require a time series of both \da{} and \id{} prices for which, to the best of our knowledge, no forecasts exists.
Therefore, we pragmatically consider the time series of the years 2020, 2021, and 2022 in our analysis, assuming that these represent scenarios for  low, medium, and high future electricity prices.

To compute the corresponding historical wind and PV power time series, we use weather data for Aachen, Germany, from the German weather service \citep{DeutscherWetterdienst,wetterdienst_python}.
Specifically, we pre-process the measured wind speed, global radiation, and diffuse radiation similar to \cite{Bahl2017} and obtain the relative PV and wind power generation as discussed in detail in Section 2 of the supporting material. 
Section 1 of the supporting material shows the mean and the standard deviation of the historical wind speed and solar irradiance that stay within rather narrow ranges, with 2020 as a windier year and 2022 as a sunnier year. 

Following \cite{Baumgartner2019} and \cite{Nilges2023}, we determine the average emission factor of electricity from the German grid for every time step, i.e., $\text{GWI}^\text{el}_{s,t}$ used in \cref{eqn:GWI}, by considering the momentary mix of power sources and their respective emission factors based on data of \cite{SMARD} and the ecoinvent database \citep{Ecoinvent}, respectively. 

\begin{figure}[!ht]
\centering
  \includegraphics[trim=0cm 13.5cm 12.5cm 0cm ,clip, width =0.45\textwidth]{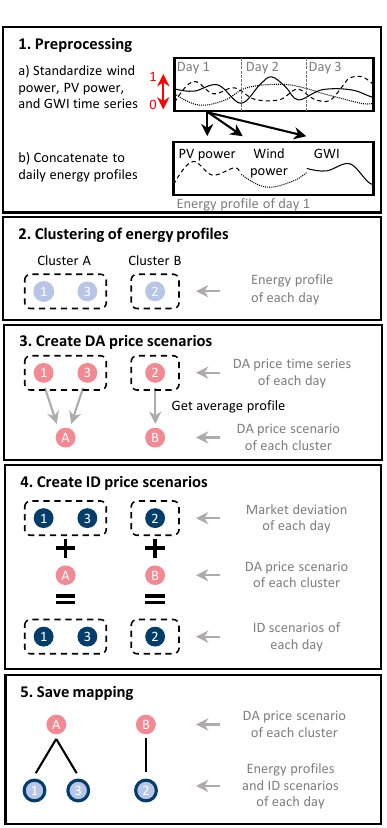}
\caption{Approach to determine scenarios based on historical time series data: The colored dots refer to the multi-dimensional data points that represent (concatenated) time series data.
\label{fig:clustering}
} 
\end{figure}

We determine the scenarios and the mapping between second and third stage of the optimization problem based on a clustering as depicted in \cref{fig:clustering}.
First, we preprocess the data by standardizing the wind and PV power time series and the GWI time series using the z-score \citep{Mohamad2013} and concatenating the daily profiles of wind power, PV power, and GWI time series to daily energy profiles.
Next, we apply k-means clustering  using the scikit-learn module in Python \citep{scikit-learn} treating each daily energy profile as a multi-dimensional data point. This leads to a clustering of similar energy profiles.
We transfer the obtained clustering to the \da{} price and market deviation time series and create an average \da{} profile as DA price scenario for each cluster.
We add the market deviations from all constituents of a respective cluster to the average \da{} price scenario to obtain \id{} price scenarios. Note that using the market deviation allows accounting for the inter-market correlation similar to our prior work \citep{Germscheid_2022}.
Finally, we use the mapping resulting from the clustering to connect the second and third stage of the optimization problem. 
Note that on the second stage, the probability of the \da{} realizations depends on the cluster size.
In contrast, the realizations on the third stage are equi-probable,  i.e., $\pi_s=1/|\mathbb{S}|$, as the clustering is used on the third stage to establish the mapping but not for data reduction.
We show in Section 3 of the supporting material that the within-cluster sum-of-squares does not allow deriving an obvious decision on a suitable number of clusters for the given data.
For our application, we look for a compromise between the number of clusters and number of scenarios per cluster.
Pragmatically, we consider 20 clusters, which leads to roughly 18 scenarios per cluster on average. We will discuss the impact of clustering on the results in \cref{sec:simultanousmarket}.

\subsection{Model specifications and evaluation}\label{sec:specifications}

\cref{tab:DoF} specifies the degrees of freedom of the integrated design and scheduling problem.
\cref{tab:process} lists reference parameters of the generalized process that we also used in our prior work \citep{Germscheid_2022}. Note that the reference parameters are similar to the chlor-alkali electrolysis \citep{Germscheid_2022} with the exception of stricter ramping limitations of the considered reference process. 
This stricter limitation allows analyzing the impact of ramping  restrictions on the potential savings of the integrated design and scheduling in the parameter study in \cref{sec:synergy}.

\begin{table}[ht!]
\caption{Degrees of freedom: 
The number of degrees of freedom depends on the number of clusters $n_c$, the number of quarter-hourly time steps $T$, and the number of scenarios $|\mathbb{S}|$, with $T=96$ and $|\mathbb{S}|=365$ in our case.
The operational degrees of freedom concern the charging and discharging of the battery.
Note that other optimization variables, e.g., power intake, PV power, and wind power, are not degrees of freedom, but can be determined from equality constraints.
}
\label{tab:DoF}
\centering
\begin{tabular}{lrr}
\hline
                                        & Degrees of freedom & Explanation \\
\hline
Design & 3      &Capacities of PV, wind power, and battery            \\
DA market  & $n_c\cdot T/4 $    & DA electricity purchases and sales     \\
ID  market  & $|\mathbb{S}| \cdot T $   &   ID electricity purchases and sales   \\
Operation  & $2\cdot |\mathbb{S}| \cdot T  $ &Battery charging and discharging rate      \\
\hline
\end{tabular}
\end{table}

\begin{table}[!ht]
    \caption{Parameters of the generalized process model: Process oversizing, minimal part load, and ramping limit are given related to the nominal power intake. The product storage capacity refers to the time necessary to fill an empty product storage considering production at nominal power intake.    }
    \centering
    \begin{tabular}{llll}
    \hline
     Parameter &  Reference values \\
     \hline
      Nominal power intake &  2.74 MW\\
      Process oversizing  &  20\% \\
      Minimal part load  &  50\% \\
      Product storage capacity   & 3h\\
      Ramping limit  & 25\%/h\\
     \hline
\end{tabular}
    \label{tab:process}
\end{table}

We consider the parameters given in \cref{tab:CAPEX} for the CAPEX.
Similar to \cite{Sass2020}, we consider an interest rate $\gamma_{1}=8\%$. 
In Section 4 of the supporting material, we list the resulting annual PV and wind electricity generation cost, showing that wind power has lower production cost  than PV due to a higher average utilization rate.
For the \gwi{} of the electricity generation and storage system, we use data of the ecoinvent database 3.9.1 \citep{Ecoinvent} that we specify in Section 4 of the supporting material for reproducibility. 
For the battery, we consider a charging and discharging rate of 4h \citep{Tesla2023}, i.e., $\tau= 4$h,  and a round-trip efficiency of $90\%$ \citep{DENA}, i.e., $\sqrt{\eta_\text{in}}=\sqrt{\eta_\text{out}}=\sqrt{0.9}$. 
Moreover, we consider the average grid fee cost of 2022 for industrial consumers in Germany, i.e., $c_\text{Fee} = 29.6$ EUR/MWh \citep{MonitoringBericht_BNA}.
Furthermore, we choose the nominal capacity of the power-intensive production process such that the process is classified as an industrial consumer, i.e., 24 GWh per year \citep{MonitoringBericht_BNA}, which allows the process operator to benefit from lower grid fees compared to non-industrial consumers \citep{MonitoringBericht_BNA}.

\begin{table}[ht!]
\caption{Component-specific life time, investment and maintenance cost for Germany based on \cite{DENA}.\label{tab:CAPEX}}
\centering
\begin{tabular}{lrrr}
\hline
    & Lifetime $\gamma_{2,i}$ & $\text{CAPEX}^0_i$ & Annual maintenance cost  $\gamma_\text{2,i}$   \\
\hline
Roof-top PV             & 25 a & 927 EUR/kWp      & 17 EUR/kWp            \\
Onshore wind         & 25 a& 1113 EUR/kWp         & 13 EUR/kWp            \\
Battery  & 15 a & 550 EUR/kWh         & 20 EUR/kWh           \\
 \hline 
\end{tabular}
\end{table}

We expect the maximum allowed capacities of the energy system to have an impact on the optimization results. 
Pragmatically, we first restrict the admissible capacities for wind power and PV by the nominal power intake, i.e., $Q_\text{W}^\text{max}=Q_\text{PV}^\text{max}=P_\text{nom}$. We choose the admissible battery size such that the maximum discharge rate corresponds to the nominal power intake of the production process, i.e.,  $P_\text{nom} = Q_\text{B}^\text{max}/ \tau$.
In \cref{sec:synergy}, we then analyze the impact of the maximum allowed energy system capacities on the TAC in detail.

We implement the model in Pyomo \citep{Pyomo} and use the solver Gurobi 9.5.0 \citep{GurobiOptimization2020} with default settings on an Intel Core i7-9700 processor and 32GB RAM. We formulate the multi-stage problem by means of its deterministic equivalent.

\section{Results}\label{sec:results}

In the following, we analyze the synergies between the local energy system and the flexible production process and the benefit of considering simultaneous market participation at design stage.
To this end, we first consider market participation only in the \id{} market and discuss the optimal design and savings of the local energy system (\cref{sec:basicDesign}) as well as the impact of the process flexibility on the potential savings (\cref{sec:synergy}).
Note that we select the \id{} market instead of the \da{} market, as the \id{} market allows adapting the electricity procurement in response to realization of the uncertainty in the renewable electricity supply.
We then show the difference between single and simultaneous market participation and discuss the energy system design in the context of \dr{} scheduling with simultaneous \da{} and \id{} market participation (\cref{sec:simultanousmarket}).

\subsection{Design and operation for single market participation} \label{sec:basicDesign}

In the following, we evaluate the energy system design considering only the \id{} market for the reference process defined in \cref{tab:process}.

\begin{figure}[!ht]
         \centering
         \includegraphics[width=\textwidth]{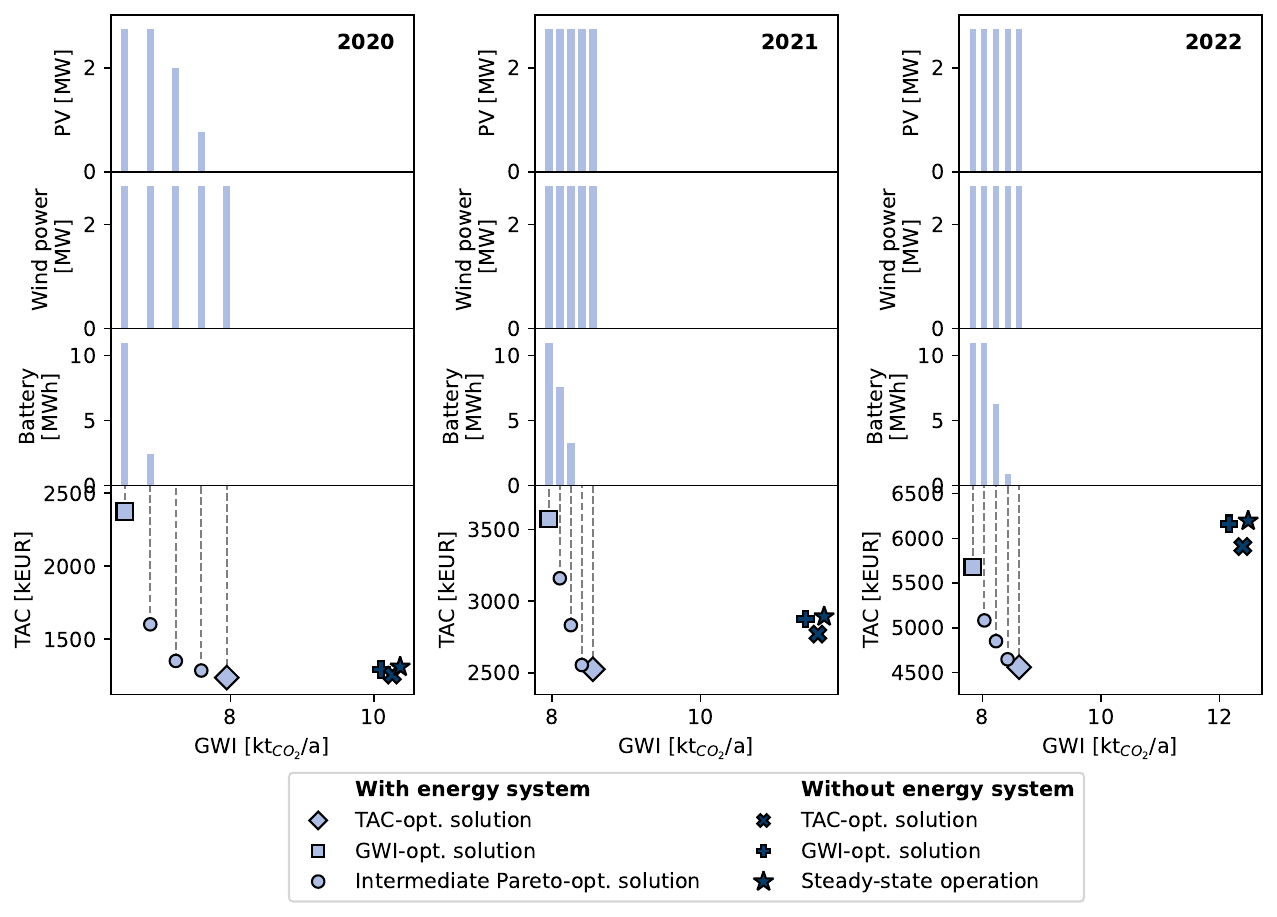}
         \caption{Energy system design for \id{} market-only participation: Pareto-optimal solutions are given for 2020 (left), 2021 (center), and 2022 (right). For each year, the \tac{}-optimal, the \gwi{}-optimal, and three intermediate Pareto-optimal solutions are shown, which are equi-distant with respect to \gwi{}. 
         \gwi{} and \tac{} (lower part) are given with vertical gray dashed lines pointing to the respective optimal capacities of the local electricity supply system (upper three parts).
         Additionally, the \tac{}- and \gwi-optimal \dr{}  as well as the steady-state operation without a local energy system  are given for comparison.}
         \label{fig:basicDesign}
\end{figure}

\cref{fig:basicDesign} shows Pareto-optimal energy system designs based on the time series for 2020, 2021, and 2022.
In three cases, the ecologic and economic objectives lead to competing solutions.  Interestingly, the \tac{}-optimal solutions do not contain a battery, as potential savings from operating the battery do not outweigh the battery investment cost.
High electricity prices in 2021 and 2022 incentivize both on-site wind and PV electricity generation in the \tac{}-optimal solutions, whereas only wind electricity generation is used in 2020. The preference for wind can be explained by the higher average utilization rate (see Section 4 of the supporting material). 
In contrast, battery, PV, and wind generation capacities are built at maximum capacity for all \gwi{}-optimal solutions, irrespective of the year studied. However, integrating a battery leads to a large increase in TAC with only a small improvement in GWI, as can be seen from the shape of the Pareto front.

\begin{figure}[!ht]
     \begin{subfigure}[b]{0.45\textwidth}
         \centering
         \includegraphics[width=\textwidth]{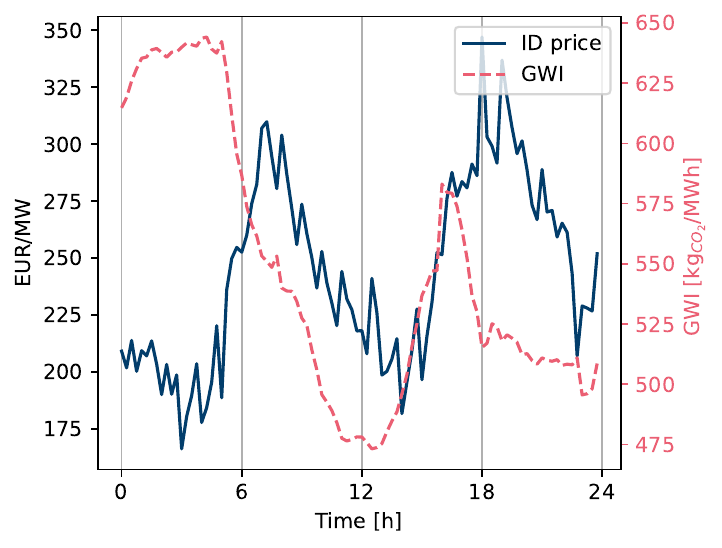}
         \caption{Electricity price and emission factor.}
         \label{fig:operationa}
     \end{subfigure}
    \hfill
     \begin{subfigure}[b]{0.45\textwidth}
         \centering
         \includegraphics[width=\textwidth]{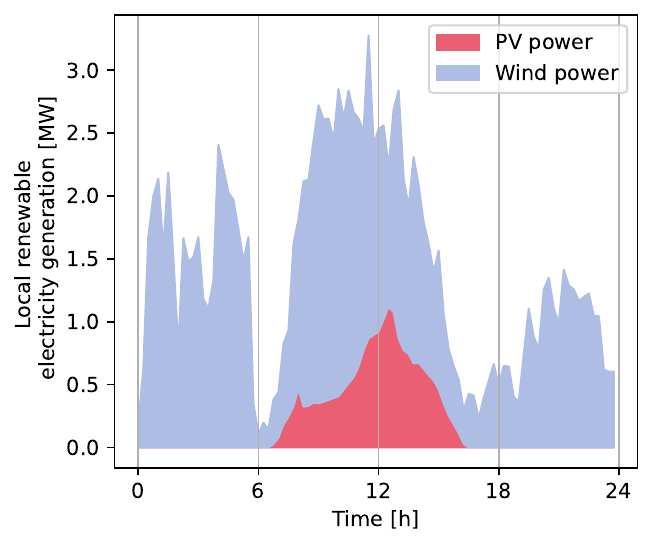}
         \caption{Renewable electricity generation.}
         \label{fig:operationb}
     \end{subfigure}

     \begin{subfigure}[b]{0.45\textwidth}
         \centering
         \includegraphics[width=\textwidth]{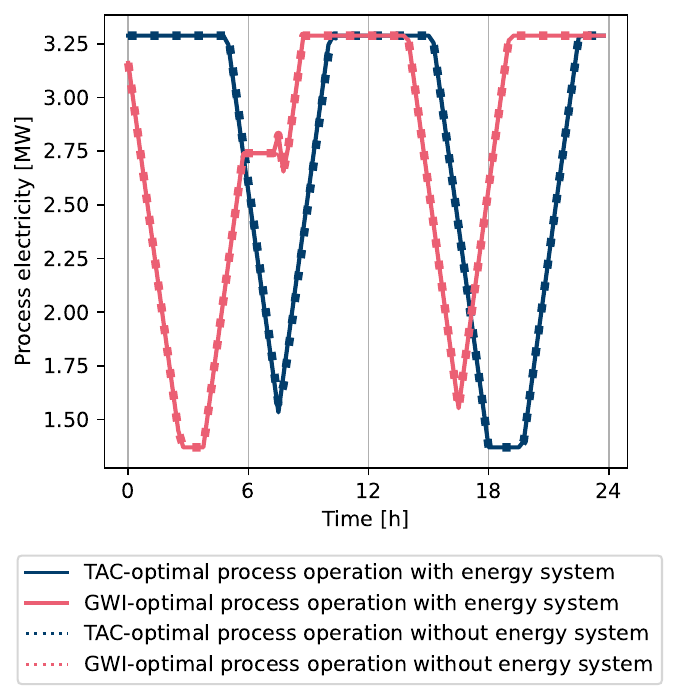}
         \caption{Operation of the flexible reference process.}
         \label{fig:operationc}
     \end{subfigure}
     \hfill
     \begin{subfigure}[b]{0.45\textwidth}
         \centering
         \includegraphics[width=\textwidth]{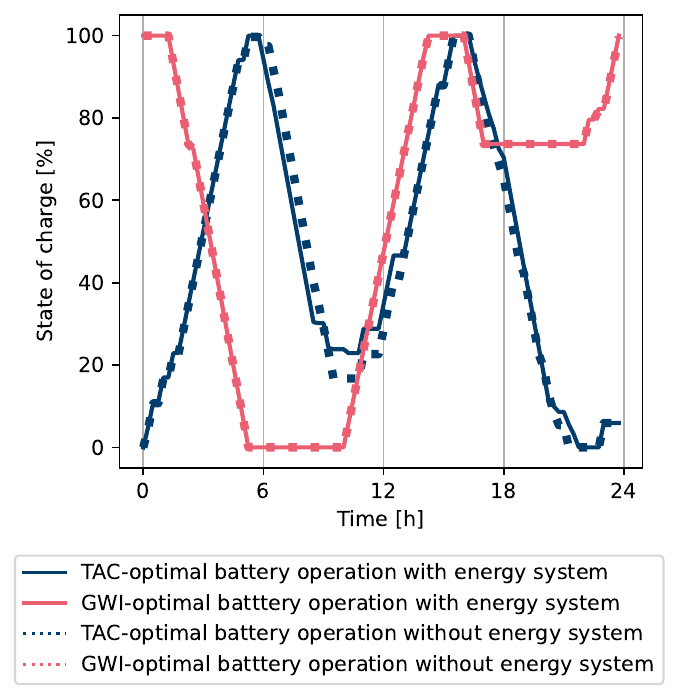}
         \caption{Battery operation.}
         \label{fig:operationd}
     \end{subfigure}  
     \caption{Process and energy system operation for an exemplary day in 2022 (\id{}-only market participation): 
     An energy system design with maximum admissible PV, wind power, and battery capacities is considered. 
     The operation is shown for the \id{} electricity price and emission factor (a) and the renewable electricity generation (b). 
    The process operation determines the electricity consumption (c). 
    The battery operation determines the state of charge (d). 
     }
     \label{fig:operation}
\end{figure}

\cref{fig:basicDesign} reveals that for 2022, the two designs with the lowest GWI have identical PV, wind power, and battery capacities but have notable differences in TAC and GWI, indicating different operating strategies.
\cref{fig:operation} shows the operation for an exemplary day given a fixed energy system design and confirms these findings. In particular, the electricity price and grid emission factor (\cref{fig:operationa}) set different incentives for TAC-optimal and GWI-optimal operation of the process (\cref{fig:operationc}) and the battery (\cref{fig:operationd}).
Note that an alignment of electricity price and grid emission factor could be achieved by increasing renewable energy penetration in the grid as well as a sufficiently high CO$_2$ price \citep{Nilges2024}.   
In case of a high renewable energy penetration in the grid, the benefit of on-site renewable generation would decrease. However, on-site generation would still be advantageous as it partially avoids grid fee cost and electricity losses due to long-distance transmission, reduces the need to expand the electricity grid, and helps avoiding the emissions and costs associated to grid expansion.

\cref{fig:operationc} additionally shows the process operation without a local energy system, revealing a similar \dr{} schedule as for a process with a local system.
Similarly, \cref{fig:operationd} shows the battery operation with and without local renewable electricity generation revealing a similar operating pattern with only minor differences.
We attribute this behavior to the time-varying incentives for \dr{} at the operational level, i.e., the electricity price and the grid emission factor. The incentives predominantly influence the operation of the process and the local energy system, while on-site generated electricity has a minor influence.

\cref{fig:basicDesign} suggests that the difference between steady-state operation and \dr{} without local energy system remains rather similar. In contrast, the difference between \dr{} without energy system and \dr{} with a local energy system  increases each year in particular with respect to the TAC. 
\cref{tab:ID_relsaves} compares ecologic and economic savings, i.e., savings with respect to the GWI and cost, respectively,  resulting from \dr{} and the local energy system and confirms these findings.
The absolute economic savings from \dr{} in comparison to steady-state operation increase due to the increasing standard deviation of the electricity price (\cref{fig:ID_price}).
The relative and absolute economic savings resulting from the energy system increase due to the increasing grid electricity cost. In particular, the savings increase significantly from 1.4\% in 2020 to 22.8\% in 2022. 
Looking at the GWI-optimal solution, the ecologic savings resulting from the energy system are much larger than the savings from \dr{} in comparison to steady-state operation. 
Furthermore, the ecologic savings are somewhat similar in all years, i.e., between 30.3\% and 35.6\%, and the variance can be attributed to the natural variability of PV and wind power production and the varying grid emission factor.

\begin{table}[ht!]
\caption{Economic and ecologic savings due to \dr{} and local electricity generation and storage.
The savings are compared for the \tac{}-optimal solutions (upper part) and the \gwi{}-optimal solutions (lower part).
Single market participation in the \id{} market is considered.
}
\centering
\begin{tabular}{lrrrr} \hline
 & 2020 & 2021 & 2022\\\hline
\multicolumn{2}{l}{\textbf{TAC-optimal solution} } &  \\
 \hline  \multicolumn{4}{l}{DR vs steady-state operation  (no energy system)} &  \\
Relative economic savings & 4.3\% & 4.2\% & 4.6\% \\
Absolute economic savings [kEUR] & 56 & 123 & 288 \\
\hline \multicolumn{3}{l}{DR with energy system vs DR without energy system } &  \\
Relative economic savings & 1.4\% & 8.8\% & 22.8\% \\
Absolute economic savings [kEUR] & 17 & 245 & 1349 \\
\hline \hline  \multicolumn{2}{l}{\textbf{GWI-optimal solution} } &  \\
\hline \multicolumn{4}{l}{DR vs steady-state operation  (no energy system)} &  \\
Relative ecologic savings & 2.5\% & 2.2\% & 2.6\% \\
Absolute ecologic savings [kt$_{{CO}_2}$/a] & 0.3 & 0.3 & 0.3 \\
\hline \multicolumn{3}{l}{DR with energy system vs DR without energy system } &  \\
Relative ecologic savings & 35.4\% & 30.3\% & 35.6\% \\
Absolute ecologic savings [kt$_{{CO}_2}$/a] & 3.6 & 3.5 & 4.3 \\
\hline
\end{tabular}
\label{tab:ID_relsaves}
\end{table}

\FloatBarrier

\subsection{Parameter study of process flexibility and energy system capacity} \label{sec:synergy}

In the following, we analyze the impact of process flexibility and admissible energy system capacities on TAC and GWI. 

\begin{figure}[!ht]
     \begin{subfigure}[b]{0.45\textwidth}
         \centering
         \includegraphics[width=\textwidth, valign=t]{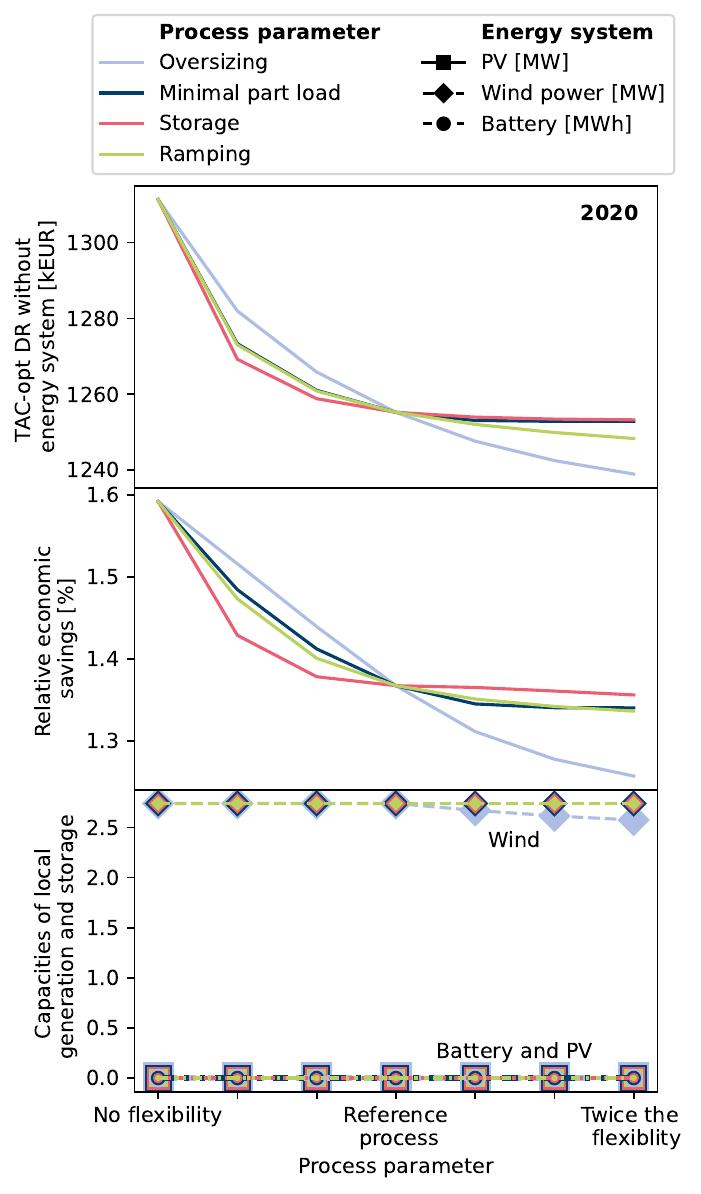}
         \caption{TAC-optimal solutions for 2020.}
         \label{fig:synergya}
     \end{subfigure}
     \hfill
     \begin{subfigure}[b]{0.45\textwidth}
         \centering
         \includegraphics[width=\textwidth, valign=t]{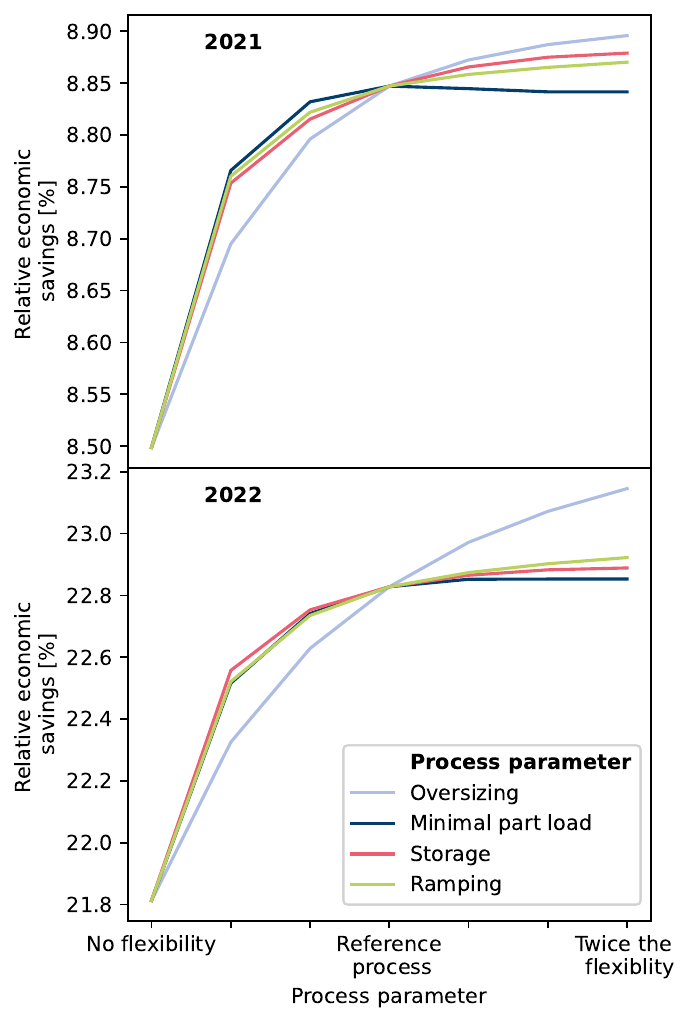}
         \caption{Relative economic savings due to a local energy system for 2021 and 2022}
         \label{fig:synergyb}
     \end{subfigure}
     \caption{Economic impact of different degrees of process flexibility: The TAC-optimal cost of \dr{} without energy system (a,\;top), the savings due to a local energy system (a,\;center), and optimal energy system design (a,\;bottom) for 2020 are shown. Furthermore, the savings due to a local energy system are shown for 2021 (b,\;top) and 2022 (b,\;bottom). 
     }
     \label{fig:synergy}
\end{figure}

For a parameter study on the process flexibility, we fix all process parameters to their respective reference values as defined in \cref{tab:process} and vary one process parameter at a time between the value corresponding to an inflexible process and the value corresponding to twice the flexibility of the process parameter.
\cref{fig:synergy} shows the impact of varying degrees of process flexibility on the \tac{}-optimal \dr{} without a local energy system, the economic savings enabled by a local energy system, and the optimal capacities of the local electricity generation and storage. \cref{fig:synergya} (top) shows exemplary for 2020 that without a local energy system, the flexible process particularly benefits from oversizing. 
Behaviors for 2021 and 2022 are similar and thus the corresponding figures are omitted.
The results confirm the findings of our prior work \citep{Germscheid_2022}, where we considered price data of 2019. 
Furthermore, \cref{fig:synergya} (bottom) shows the optimal capacities of the energy system components and reveals that the optimal wind capacity for 2020 slightly decreases with a high degree of process oversizing. Thus, process flexibility can impact the optimal energy system capacity. 
Corresponding figures for 2021 and 2022 are omitted, as no impact of the process flexibility on the resulting optimal designs can be found.

\cref{fig:synergy} reveals that the range of the economic savings due to a local energy system for any given year is rather narrow, i.e., varying the process flexibility does not impact the relative savings in a strong manner. 
Even though the range is small, varying the process oversizing has the largest leverage on the savings in comparison to the other process parameters.
Interestingly, the impact of the process parameters may differ depending on the investigated year, as higher process flexibility actually leads to lower relative savings in 2020 (\cref{fig:synergya}, center) whereas higher relative savings are recorded for 2021 and 2022 (\cref{fig:synergyb}). 
However, an analysis of the cost contributions shows exemplary for varying oversizing that the absolute TAC monotonously decreases, irrespective of the case with or without a local energy system (see Section 5 of the supporting material).

\begin{figure}[!ht]
    \centering
    \includegraphics[width=\textwidth]{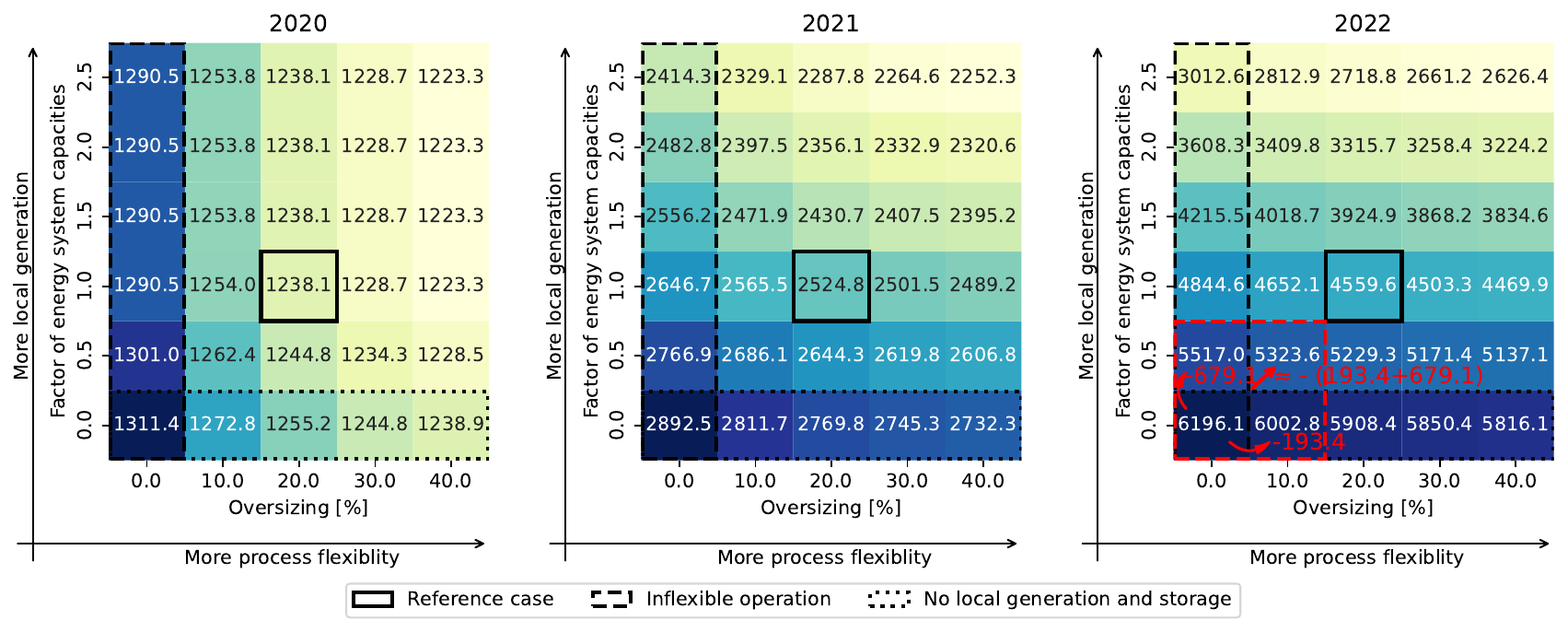}  
     \caption{Optimal TAC [EUR] with varying process oversizing and maximum admissible energy system capacity for 2020 (left), 2021 (center), and 2022 (right):  For 2022, the approximately additive behavior of the absolute savings is shown exemplary in red.
     }
     \label{fig:heatmap_TAC}
\end{figure}

\cref{fig:heatmap_TAC} shows the optimal TAC for varying admissible energy system size and process oversizing, the latter having the largest flexibility leverage for the economic savings. 
Specifically, we vary the maximum allowed energy system capacities, i.e., $Q_\text{W}^\text{max}$, $Q_\text{PV}^\text{max}$, and $Q_\text{B}^\text{max}$,  by a joint scaling factor. 
Corresponding optimal capacities of the energy system can be found in Section 5 of the supporting material. 
For 2020, process oversizing has a larger leverage on the TAC than local electricity generation and storage. 
Furthermore, the TAC remains constant for scaling factors larger than one, as a cost-optimal maximum of the wind power capacity is attained (see Section 5 of the supporting material). 
For 2021 and 2022, it can be seen that local electricity generation and storage is more economically attractive than process flexibility.

Interestingly, the absolute savings of the TAC-optimal solution enabled by a higher process flexibility and by a larger energy system behave roughly additively, which is  shown exemplary for 2022 in \cref{fig:heatmap_TAC}.
Section 5 of the supporting material shows the relative difference between the optimal TAC and the estimated TAC, the latter being defined as the sum of the absolute savings from process flexiblization and installation of a local energy system. The difference being rather small, i.e., less than 0.5\%, means that a quick first  approximate economic assessment can be performed by considering the savings from \dr{} and the cost savings from a local energy system independently.

\begin{figure}[!ht]
\centering
     \begin{subfigure}[b]{0.45\textwidth}
         \includegraphics[width=\textwidth]{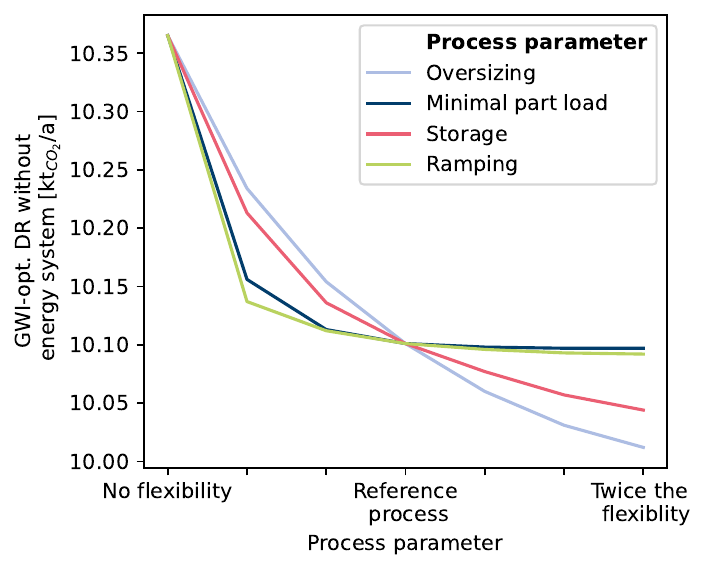}
         \caption{GWI-optimal \dr{} without a local energy system}
         \label{fig:gwi_market}
     \end{subfigure}
    \hfill
     \begin{subfigure}[b]{0.45\textwidth}
         \includegraphics[width=\textwidth]{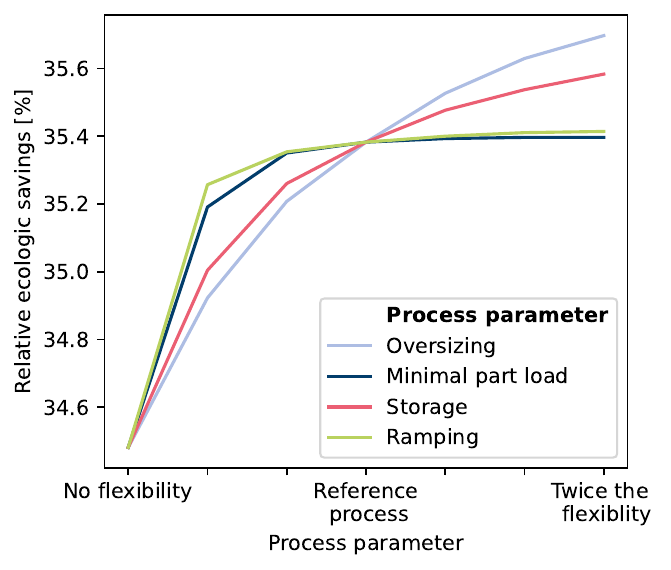}
         \caption{Relative ecologic savings due to a local energy system }
          \label{fig:gwi_savings}
     \end{subfigure}
     \caption{Ecologic impact of varying process flexibility for 2020: GWI-optimal \dr{} without a local energy system (a) and relative savings due to installation of a local system (b) are shown.  }
     \label{fig:synergy_gwi}
\end{figure}

\cref{fig:synergy_gwi} shows the impact of varying degrees of process flexibility on the \gwi{}-optimal solution exemplary for 2020. Figures for the other years show similar behavior and are therefore omitted.
\cref{fig:gwi_market} shows the case of \dr{} without a local energy system and reveals that the process oversizing and the product storage capacity have the largest impact on the GWI. This finding is consistent with the results of our prior work \citep{Schafer2020_generic}, where we considered the residual load as ecologic objective instead of the GWI. 
\cref{fig:gwi_savings} shows that process oversizing and product storage capacity have the largest leverage on the ecologic savings. 
Note that, in general, the degree of process flexibility has a rather low impact on the range of the ecologic savings.

\FloatBarrier

\subsection{Simultaneous market participation}\label{sec:simultanousmarket}

Finally, we evaluate the benefit of considering  simultaneous \da{} and \id{} market participation in the integrated design and scheduling problem. \cref{tab:runtime} compares the average wall-clock run times of the single and the simultaneous market participation. As expected, the single market participation is slightly faster than the simultaneous participation that additionally contains the \da{} trading decisions.

\cref{tab:purchases} lists the TAC-optimal designs for the reference process.
Designs for 2020 and 2021 do not reveal differences to the case of \id{}-only market participation, i.e., the price differences between the markets in these years are not large enough to incentivize the installation of a battery. 
In contrast, a battery is built for the simultaneous market participation in 2022 as the battery offers trading capacity for exploitation of large price differences between the \da{} and \id{} market. 
The Pareto-optimal energy system designs of 2022 vary with respect to the battery capacity for the simultaneous \da{} and \id{} market participation compared to the \id{}-only case and are shown in Section 6 of the supporting material.

\begin{table}[ht!]
\caption{Run time comparison of single and simultaneous market participation: The stated wall-clock times of the integrated design and scheduling problem are the averages over 15 Pareto-optimal solutions (5 Pareto-optimal solutions per year).
The solver Gurobi 9.5.0 \citep{GurobiOptimization2020} was used and the machine was equipped with an Intel Core i7-9700 processor and 32GB RAM.
}
\label{tab:runtime}
\centering
\begin{tabular}{lr}
\hline
    & Run time       \\
\hline
Single market participation &        21.8 s   \\
Simultaneous market participation &   29.6 s    \\ 
\hline
\end{tabular}
\end{table}

\begin{table}[ht!]
\caption{System capacities and electricity trading amounts for ID-only participation  (upper part) and simultaneous \da{} and \id{} market participation (lower part): In all cases, TAC-optimal  \dr{} of the reference process (\cref{tab:process}) with a local energy system is considered.
The stated values of the \id{}-only participation correspond to the \tac{}-optimal solution from \cref{fig:basicDesign}.
}
\centering
\begin{tabular}{lrrr} \hline
 & 2020 & 2021 & 2022\\
 \hline  \multicolumn{3}{l}{\textbf{ID-only participation} } &  \\
\hline PV capacity & - & 2.74 MW & 2.74 MW \\
Wind capacity & 2.74 MW & 2.74 MW & 2.74 MW \\
Battery capacity & - & - & - \\
\hline Total purchases & 17,316 MWh & 16,008 MWh & 14,776 MWh \\
Total sales & 62 MWh & 203 MWh & 357 MWh \\
\hline \hline  \multicolumn{3}{l}{\textbf{Simultaneous market participation} } &  \\
\hline PV capacity & - & 2.74 MW & 2.74 MW \\
Wind capacity & 2.74 MW & 2.74 MW & 2.74 MW \\
Battery capacity & - & - & 5.31 MWh \\
\hline DA purchases & 18,582 MWh & 16,389 MWh & 18,874 MWh \\
ID purchases & 6,796 MWh & 7,442 MWh & 12,944 MWh \\
Total purchases & 25,378 MWh & 23,831 MWh & 31,818 MWh \\
\hline DA sales & 452 MWh & 954 MWh & 5,947 MWh \\
ID sales & 7,678 MWh & 7,080 MWh & 11,083 MWh \\
Total sales & 8,130 MWh & 8,034 MWh & 17,030 MWh \\
\hline
\end{tabular}
\label{tab:purchases}
\end{table}

\cref{tab:purchases} compares single and simultaneous market participation with respect to the electricity purchases and sales. 
It shows that for \id{}-only participation more electricity is purchased  in 2020, compared to 2021 and 2022, which are years with larger optimal PV and wind power capacities.
For the simultaneous participation in 2020, the majority of the electricity is purchased on the \da{} market due to the positive price difference (\cref{fig:Dev_price}).
The combination of PV and wind enables higher \da{} sales in 2021.
Moreover, total purchases and sales significantly increase for the simultaneous participation in 2022 due to an increased trading capacity enabled by the battery.

\cref{tab:DAID_ID} shows that the relative savings of simultaneous market participation compared to single market participation stay within a similar range, irrespective of the considered year.
In contrast, the absolute savings increase each year due to the increased variance of the market deviation (\cref{fig:Dev_price}). 
\cref{tab:DAID_ID} reveals that both the absolute and relative savings of simultaneous market participation increase with both process flexibilization and the integration of local electricity generation and storage.

\begin{table}[ht!]
\centering
\caption{
The TAC savings from simultaneous market participation compared to single market participation for an inflexible process without an energy system (top), a flexible process without an energy system (center), and a flexible process with an energy system (bottom). The flexible process is the reference process (\cref{tab:process}).
}
\begin{tabular}{lrrr} \hline
 & 2020 & 2021 & 2022\\\hline
 \multicolumn{4}{l}{\textbf{Inflexible process without energy system} }  \\
Relative savings & 3.0\% & 2.0\% & 2.0\% \\
Absolute savings & 39.3 kEUR & 57.4 kEUR & 122.7 kEUR \\
\hline
 \multicolumn{4}{l}{\textbf{Flexible process without energy system} }  \\
Relative savings & 3.8\% & 2.5\% & 2.5\% \\
Absolute savings & 47.1 kEUR & 68.9 kEUR & 147.2 kEUR \\
\hline
  \multicolumn{4}{l}{\textbf{Flexible process with energy system} } \\
Relative savings & 3.9\% & 2.9\% & 4.4\% \\
Absolute savings & 47.9 kEUR & 72.8 kEUR & 198.7 kEUR \\
\hline
\end{tabular}
\label{tab:DAID_ID}
\end{table}

\FloatBarrier

\cref{tab:savings_DAID} shows the contributions of cost savings for the year 2022 considering TAC-optimal simultaneous market participation. Here, an inflexible process without an energy system is modified by separately adding a battery, renewable electricity generation, and process flexiblization. 
\cref{tab:savings_DAID} reveals that the main cost savings result from the on-site electricity generation followed by process flexiblization. 
The integration of a battery accounts only for a small fraction of the savings.
Note that similar to \id{}-only participation (\cref{sec:synergy}), summing up the absolute savings from battery installation, renewable electricity generation, and flexiblization, separately, allows for a good overall savings estimate.

\begin{table}[ht!]
\centering
\caption{Savings contributions of simultaneous market participation in the TAC-optimal case for 2022: 
The savings are related to the inflexible process without an energy system (first row). 
The battery capacity of the inflexible process with the battery (third row) is identical to the  flexible process with an energy system (second row).
}
\begin{tabular}{lrrr}
\hline & TAC & Savings\\
 \hline Inflexible process  without energy system & 6073.5 kEUR & -\\
 Flexible process with energy system & 4360.9 kEUR & 1712.6 kEUR\\
 \hline Inflexible process with battery & 6044.7 kEUR & 28.8 kEUR\\
Inflexible process with wind and PV power & 4702.5 kEUR & 1371.0 kEUR\\
Flexible process  without energy system & 5761.2 kEUR & 312.3 kEUR\\
\hline\end{tabular}
\label{tab:savings_DAID}
\end{table}

\FloatBarrier

\begin{figure}[!ht]
\centering
     \begin{subfigure}[b]{0.45\textwidth}
         \includegraphics[width=\textwidth]{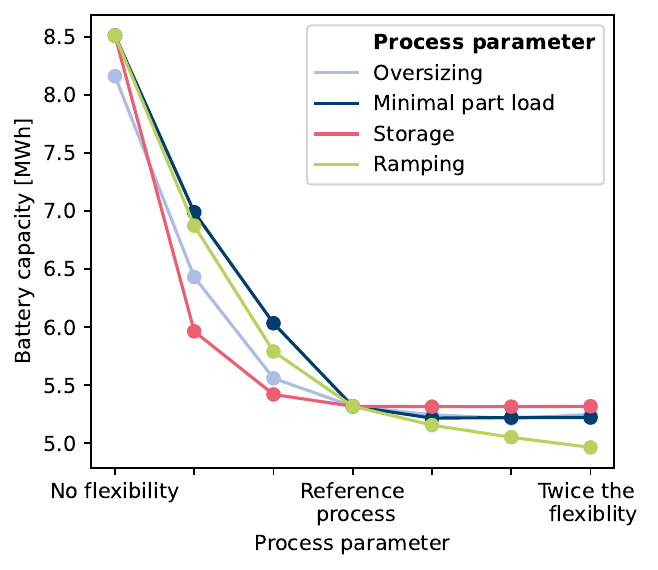}
     \end{subfigure}
\caption{TAC-optimal battery capacity for simultaneous market participation for 2022 with varying degree of process flexibility.}
\label{fig:DAID_bat}
\end{figure}

Similar to \cref{sec:synergy}, we vary the degree of process flexibility.
\cref{fig:DAID_bat} shows the impact of process flexibility on the optimal battery capacity for 2022. 
An analysis of the optimal capacities of PV and wind power is omitted as these quantities do not vary in response to a varying degree of process flexibility.
It can be noted that larger batteries are built for less flexible processes, in particular for processes with part load and ramping restrictions.
\dr{} without energy system and economic savings of \dr{} with an energy system behave similarly to the case of single market participation (\cref{sec:synergy}), irrespective of the investigated year. Thus, respective figures are omitted.

Recall that we pragmatically consider 20 clusters for the assessment. In Section 6 of the supporting material, we show the results for 10 and 30 clusters. More clusters lead to higher savings enabled by the simultaneous market participation and a larger optimal battery capacity in 2022.
In the supporting material, we show that the average standard deviation of the market deviation decreases with increasing number of clusters, which enables better positioning on the two markets in case of the simultaneous participation and, thus, increases savings and incentivizes a larger  battery.

\FloatBarrier

\section{Conclusion}\label{sec:conclusion}

We assessed the optimal design of a local electricity generation and storage system for a generalized continuous, power-intensive  production process that is capable of performing demand response and can act on both the day-ahead and intraday electricity market.
In a three-stage stochastic problem, we optimized the capacities of photovoltaic power, wind power, and electric battery  with an integrated  demand response scheduling of the production process.
Building on our prior work \citep{Schafer2020_generic,Germscheid_2022}, we used a generalized process model with few flexibility-defining parameters, i.e., process oversizing, minimal part load, product storage, and ramping limitation.
In a bi-objective optimization, we considered both economic and ecological objectives.
We considered scenarios of low, intermediate, and high electricity prices for a plant location in Germany as well as a time-varying grid emission factor.

We find that batteries are mainly built  to lower the global warming impact, however, leading to a significant increase in total annualized cost. 
Economically and ecologically-optimal operation of the process and battery primarily respond to the time-varying electricity price and grid emission factor, but only to a little extent to the on-site generation of renewable electricity.
Varying the degree of process flexibility, we find a rather small impact on the achievable relative economic and ecologic savings that come with local electricity generation and storage.
Moreover, we show that the absolute cost savings from flexiblizing the process and installing a local energy system are approximately additive. 
Comparing intraday-only and simultaneous day-ahead and intraday market participation, we find that the energy system designs are similar for the investigated scenarios, except when high price differences between the markets incentivize the installation of a battery. The cost-optimal battery capacity significantly depends on the available process flexibility, enables large volumes for trading on the markets, but comes with only minor economic savings.

In our assessment, we pragmatically considered time series based on historic data of three years to understand the effects of low, medium, and high electricity prices. 
To account for long-term variability of prices and the long lifetime of both the process and the energy system equipment, our approach should be extended to consider multiple years by incorporating long-term time series forecasting, e.g., 
based on the hourly day-ahead price forecasting of  \cite{Ziel2018,Gabrielli2022} that would need to be extended to also consider quarter-hourly intraday prices.
Uncertainties associated to the forecasts and financial risks must be particularly accounted for, e.g., in a risk-averse optimization similar to \cite{Xuan2021,Vieira2021}.
In order to use our approach to evaluate the potential of a local system for a specific process and location, the process characteristics must be known and respective local weather data is required.

\section*{Authorship contribution}

\textbf{Sonja H. M. Germscheid:} Conceptualization, Methodology, Software, Investigation, Data curation, Writing - original draft, Visualization.
\textbf{Benedikt Nilges:} Data curation, Writing – review \& editing.
\textbf{Niklas von der Assen:} Funding acquisition,  Writing – review \& editing.
\textbf{Alexander Mitsos:} Writing - review \& editing, Supervision, Funding acquisition.
\textbf{Manuel Dahmen:} Conceptualization, Methodology, Writing - review \& editing, Supervision, Funding acquisition.

\section*{Declaration of Competing Interest}
We have no conflict of interest.

\section*{Acknowledgements}
SG gratefully acknowledges the financial support of the Helmholtz Association of German Research Centers through program-oriented funding (POF) and the grant \textit{Uncertainty Quantification
– From Data to Reliable Knowledge (UQ)}  (grant number: ZT-I-0029). AM and MD acknowledge funding from the Helmholtz Association of German Research Centers through program-oriented funding (POF).
BN gratefully acknowledges the financial support of the Kopernikus project SynErgie (grant number 03SFK3L1-2) by the Federal Ministry of Education and Research (BMBF). 
This work was performed as part of the \textit{Helmholtz School for Data Science in Life, Earth and Energy (HDS-LEE)}.
We kindly thank  Yifan Wang  (RWTH Aachen University, Institute of Technical Thermodynamics) for providing the wind turbine performance curve developed by \cite{Bahl2017} which was used for pre-processing of the wind data.

\FloatBarrier
 \bibliographystyle{apalike}
  \bibliography{BibTex.bib}

\end{document}


\newacronym{da}{DA}{day-ahead}
\newcommand{\da}{\gls*{da}}
\newacronym{id}{ID}{intraday}
\newcommand{\id}{\gls*{id}}
\newacronym{dr}{DR}{demand response}
\newcommand{\dr}{\gls*{dr}}
\newacronym{tac}{TAC}{total annualized cost}
\newcommand{\tac}{\gls*{tac}}
\newacronym{gwi}{GWI}{global warming impact}
\newcommand{\gwi}{\gls*{gwi}}
\newcommand{\maxP}{\ensuremath{\theta_\text{max}}}
\newcommand{\Pnom}{\ensuremath{P_\text{nom}}}
\newcommand{\dt}{\ensuremath{\Delta{t}}}
\newcommand{\nS}{\ensuremath{\mathbb{S}}}

\ifx\REVIEW\undefined
\twocolumn[
\begin{@twocolumnfalse}
\fi
  \thispagestyle{firststyle}
  \begin{center}
    \begin{large}
      \textbf{\mytitle}
    \end{large} \\
    \myauthor
  \end{center}
  \vspace{0.5cm}
  \begin{footnotesize}
    \affil
  \end{footnotesize}
 
\ifx\REVIEW\undefined
\end{@twocolumnfalse}
]
\fi

\doublespacing

\section{Historical data}\label{sec:data_development}

\begin{figure}[!ht]
\captionsetup[subfigure]{justification=centering}
\centering
\begin{subfigure}[t]{.32\textwidth}
    \centering
  \includegraphics[width =\linewidth]{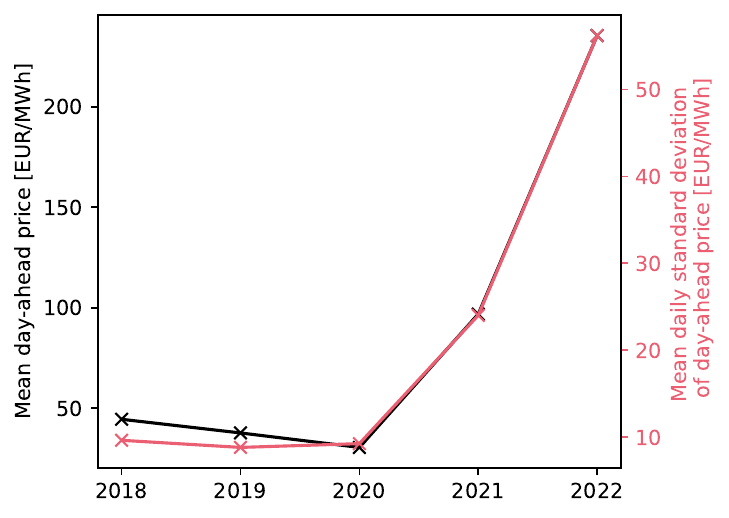}
  \caption{Day-ahead electricity price. }
  \label{fig:DA_prices}
  \end{subfigure}
  \hfill
\begin{subfigure}[t]{.32\textwidth}
    \centering
  \includegraphics[width =\linewidth]{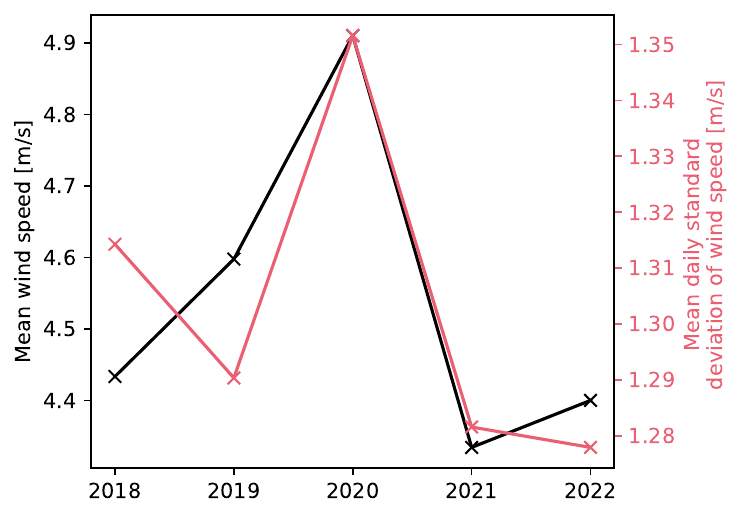}
  \caption{Wind speed measured 10m above ground.}
  \label{fig:wind}
  \end{subfigure}
\hfill
 \begin{subfigure}[t]{.32\textwidth}
    \centering
  \includegraphics[width =\linewidth]{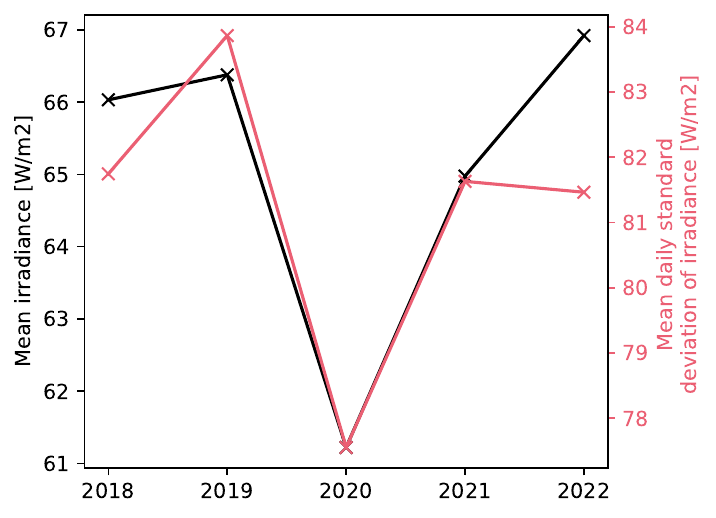}
  \caption{Solar radiance on a tilted surface.}
  \label{fig:irradiance}
  \end{subfigure}
\caption{Annual mean and mean of daily standard deviation of historical data.} 
\label{fig:data}
\end{figure}

In Section 2.4 of the main manuscript, we show the historical development of the \acrfull{id} electricity market price and the market deviation by means of the annual mean and the annual mean daily standard deviation.
\cref{fig:data} shows corresponding figures for the \acrfull{da} price (\cref{fig:DA_prices}), the wind speed (\cref{fig:wind}), and the solar irradiance on a tilted surface (\cref{fig:irradiance}). The solar irradiance on a tilted surface considers global and diffuse radiation as described in detail in \cref{sec:preprocessing} of the supporting material.
Similar to the \id{} price, the \da{} price shows a lower mean in the initial phase of the COVID-19 pandemic (year 2020) and an increased mean and standard deviation prior and during the conflict in Ukraine (year 2021 and 2022).
The standard deviation of the \da{} price is smaller than that of the \id{} price.
In comparison to the market price, relative changes of wind speed and solar irradiance are smaller.

\newpage

\section{Pre-processing of weather data}\label{sec:preprocessing}

We base the scenarios in Section 2.4 of the main manuscript on historic weather data, i.e., wind speed, global radiation, and diffuse radiation. The historic weather data have a time resolution of 10 minutes to which we apply linear interpolation to obtain a 15 minute resolution to match the assumption of constant quarter-hourly renewable electricity generation made in Section 2.1 of the main manuscript.
From the weather data, we derive the relative wind power output $\bar{q}_{\text{Wind}}$ and the relative photovoltaic (PV) power output $\bar{q}_\text{PV}$ that are necessary to calculate PV and wind power generation in Eqs. (11) and (12) of the main manuscript and are described in the following. 

\begin{figure}[!ht]
\centering
  \includegraphics[width =0.45\textwidth]{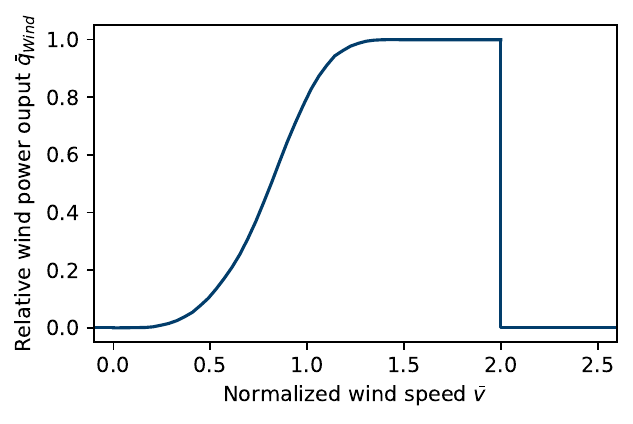}
\caption{ Performance curve of a generic wind turbine $f_\text{performance}$ by \cite{Bahl2017}: The performance curve $f_\text{performance}$ allows evaluating the relative wind power output $\bar{q}_{\text{Wind}}$ from the normalized wind speed $\bar{v}$. The curve is based on 50 data points \citep{Bahl2017}.
\label{fig:perfromance_curve}
} 
\end{figure}

We calculate the relative wind power output $\bar{q}_{\text{Wind}}$ using the performance curve of a generic wind turbine following \cite{Bahl2017}: 
\begin{align}
 v_\text{Hub}&=v_\text{Measure} \frac{\text{ln} (H_\text{Hub}/Z_0)}{\text{ln} (H_\text{Measure}/Z_0)} \label{eqn:vhub},\\
 \Bar{v}&= \frac{v_\text{Hub}}{v_\text{Ref}}\label{eqn:vrel},\\
  \bar{q}_{\text{Wind}}&= f_\text{performance}(\Bar{v}) \label{eqn:qrelqind}
\end{align}
In \cref{eqn:vhub}, the wind speed $v_\text{Hub}$ at hub height $H_\text{Hub}$ is calculated based on the measured wind speed $v^\text{Measure}_{\text{Wind}}$  at measuring height $H^\text{Measure} = 10$m \citep{DeutscherWetterdienst} and the ground roughness $Z_0$. 
We assume a hub height of $80$m.
The wind speed was recorded in an agricultural area with few houses \citep{DeutscherWetterdienst}, i.e.,  $Z_0=0.1$m \citep{ground_roughness}. 
In \cref{eqn:vrel}, the wind speed is normalized by means of the reference wind speed $v_\text{Ref} = 11.8$m/s (\cite{Bahl2017}) for the performance curve.
Finally, the generic wind turbine performance curve $f_\text{performance}$ by \cite{Bahl2017} shown in \cref{fig:perfromance_curve} is used to evaluate the relative wind power output $\bar{q}_{\text{Wind}}$.

Similar to \cite{Bahl2017}, we calculate the irradiance  $I$ on a tilted surface from the measured global and diffuse radiation using the Perez model \citep{Perez1987} implemented in the Python package pvlib \citep{pvlib}. Similar to  \cite{Sass2020}, we assume a surface tilt of $10^{\circ}$ and an azimuth of $103^{\circ}$ and calculate the relative PV power output $\bar{q}_\text{PV}$ from the irradiance $I$, efficiency $\eta=0.19$ \citep{Sass2020}, and nominal capacity $\text{cap}_\text{nom}=0.1 \; \text{kWm}^{-2}$ \citep{DENA}:
\begin{align}
    \bar{q}_\text{PV}&= \text{min} ( \frac{\eta\;I}{\text{cap}_\text{nom}},1) \label{eqn:relPVpower}
\end{align}
In \cref{eqn:relPVpower}, we limit the maximum relative power output to one when the irradiance exceeds the conversion capacity of the PV system similar to \cite{Sass2020}.

\section{Within-cluster sum-of-squares}\label{sec:clustering}

\begin{figure}[!ht]
\centering
  \includegraphics[width =0.45\textwidth]{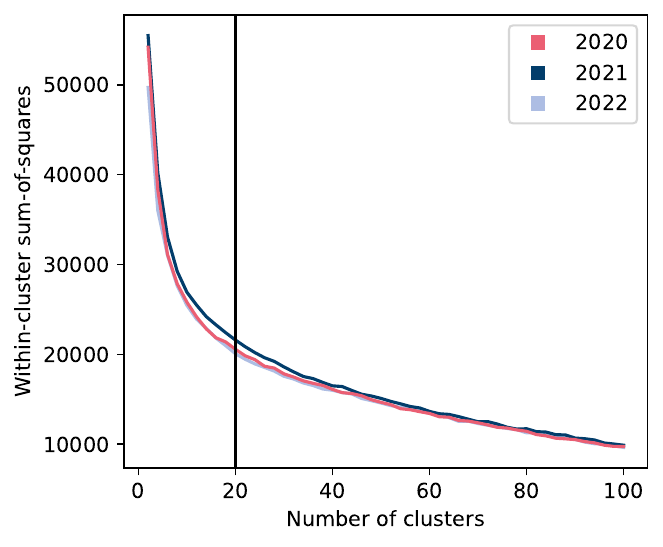}
\caption{ Within-cluster sum-of-squares are show for clustering the concatenated wind power, PV power, and grid emission time series of 2020, 2021, and 2022. 20 clusters (vertical black line) are chosen for the analysis in the main manuscript.  }
\label{fig:inertie} 
\end{figure}

In Section 2.4 of the main manuscript, we apply k-means clustering to the time series data.
\cref{fig:inertie} shows the trade-off between the within-cluster sum-of-squares, i.e., the total within-cluster variance, and the number of clusters for the concatenated wind and PV power and grid emission time series of 2020, 2021, and 2022.
The trade-off curve does not reveal a particular kink that could be used as decision criterion for a suitable number of clusters according to the elbow method \citep{Thorndike1953_elbow}.

\FloatBarrier

\section{Supporting information for the model specification}\label{sec:supp:specification}

In Section 2.5 of the main manuscript, we specify the CAPEX parameters of the local electricity generation and storage system. \cref{tab:prodPower} reports the annual PV and wind power production as well as the resulting respective electricity generation cost. \cref{tab:GWI} specifies the data used for the emission factors of the electricity generation and storage system. 
Specifications were made based on the quality of available data and the applicability to our scenarios of a plant located near Aachen, Germany.

\begin{table}[ht!]
\caption{Production cost of PV and wind power: The generated PV and wind power are given in relation to the peak capacity. The production cost consider the produced power and the annualized investment cost. The amount of produced power is subject to natural variability and thus the weather time series.}
\centering
\begin{tabular}{lrrr} \hline
 & 2020 & 2021 & 2022\\\hline
PV power  & 1.02 GWh/MWp & 1.08 GWh/MWp & 1.11 GWh/MWp \\
PV production cost & 101.6 EUR/MWh & 95.8 EUR/MWh & 93.0 EUR/MWh \\
\hline Wind power  & 2.55 GWh/MWp & 1.96 GWh/MWp & 2.4 GWh/MWp \\
Wind production cost & 46.1 EUR/MWh & 59.9 EUR/MWh & 49.0 EUR/MWh \\
\hline
\end{tabular}
\label{tab:prodPower}
\end{table}

\begin{table}[ht!]
\caption{Ecoinvent data: For the global warming impact, we use licensed data from the ecoinvent database 3.9.1 \citep{Ecoinvent}.
\label{tab:GWI}}
\footnotesize
\begin{tabular}{lll}
\hline
 & Ecoinvent specifications \\
 \hline
PV                     &   Market for photovoltaic flat-roof installation, 3kWp, single-Si, on roof (global)     \\
Wind  power        &   Market for wind turbine with network connection, 2MW turbine, onshore (global) \\
Battery                &   Market for lithium-ion battery, LiMn2O4, rechargeable, prismatic (global) \\
\hline
\end{tabular}
\end{table}

\FloatBarrier
\section{Supporting figures for the parameter study}\label{sec:supportive_figures}

\begin{figure}[!ht]
\centering
     \includegraphics[width=0.5\textwidth]{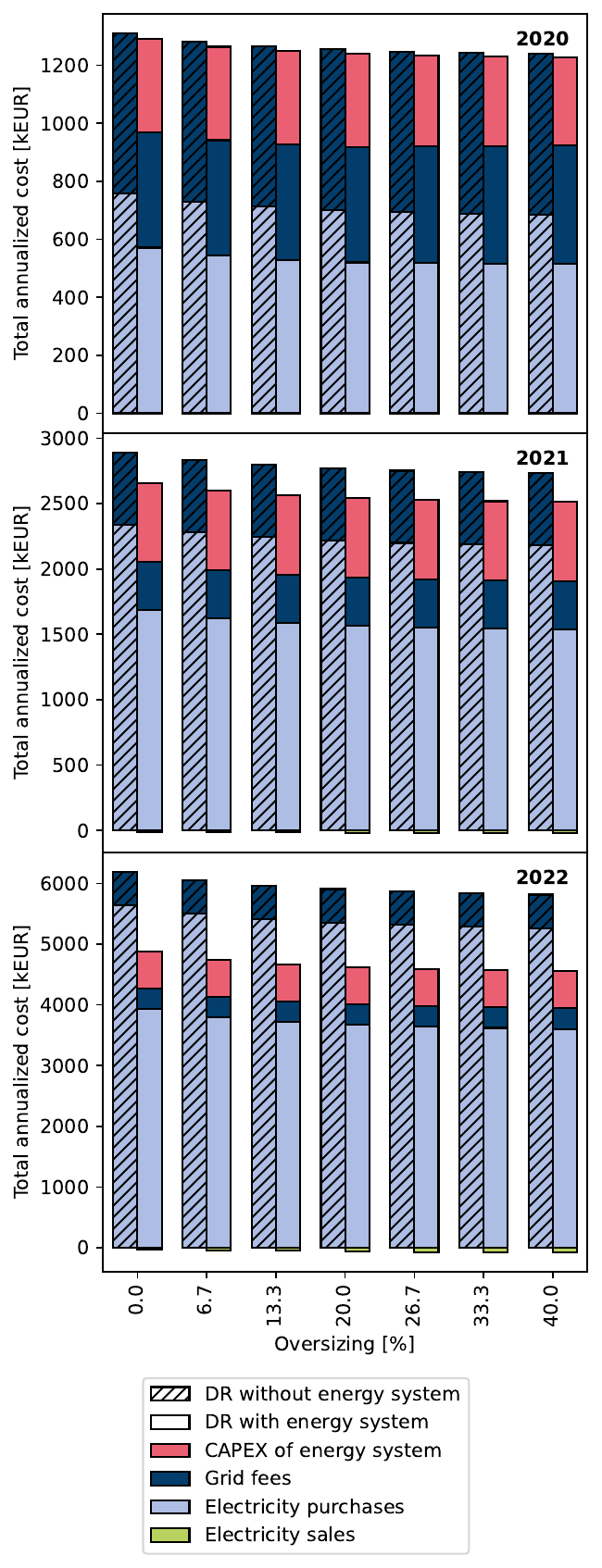}
     \caption{Cost contributions for \dr{} with energy system and \dr{} without energy system: 2020 (top), 2021 (center), and 2022 (bottom).}
     \label{fig:supp:costcontribution}
\end{figure}

Fig. 6 of the main manuscript shows the relative economic savings due to installing a local energy system for different degrees of process flexibility.
Additionally, \cref{fig:supp:costcontribution} shows the cost contributions for both optimal \dr{} without energy system and optimal \dr{} with an energy system for different degrees of process oversizing, i.e., the largest  leverage for the economic savings. 
For \dr{} without an energy system, the grid fees are constant irrespective of the considered year, while the cost of electricity purchases increases due to the higher electricity prices (Fig. 2a in the main manuscript).
Similarly, the cost of electricity purchases increases in case of \dr{} with a local energy system. Furthermore,  CAPEX of the energy system in 2021 and 2022 are larger than CAPEX in 2020 due to installation of  PV and wind power. In contrast, the  grid fee cost are lower in 2021 and 2022, compared to 2020, as less electricity is purchased (Tab. 6 in the main manuscript).
Moreover, savings due to electricity sales are higher in 2021 and 2022, compared to 2020, due to higher amount (Tab. 6 in the main manuscript) and price of sold electricity.

\begin{figure}[!ht]
     \begin{subfigure}[b]{0.32\textwidth}
         \centering
         \includegraphics[width=\textwidth]{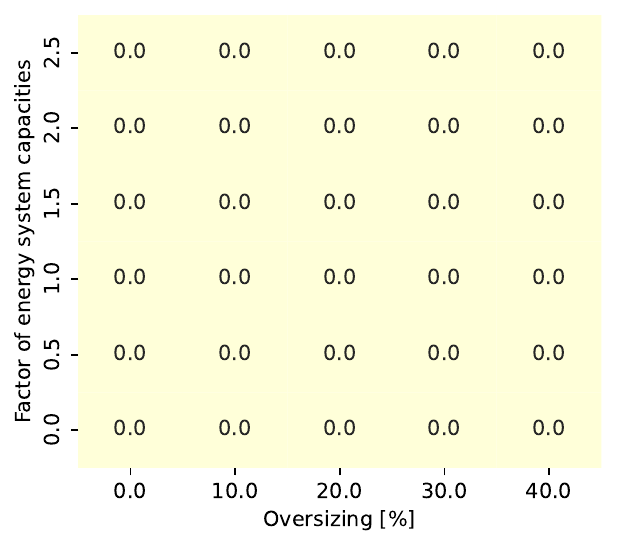}
         \caption{2020}
     \end{subfigure}
     \hfill
     \begin{subfigure}[b]{0.32\textwidth}
         \centering
         \includegraphics[width=\textwidth]{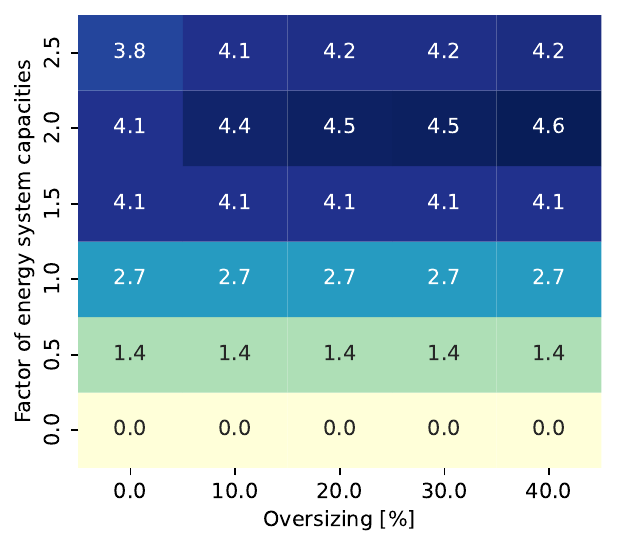}
         \caption{2021}
     \end{subfigure}
     \hfill
     \begin{subfigure}[b]{0.32\textwidth}
         \centering
         \includegraphics[width=\textwidth]{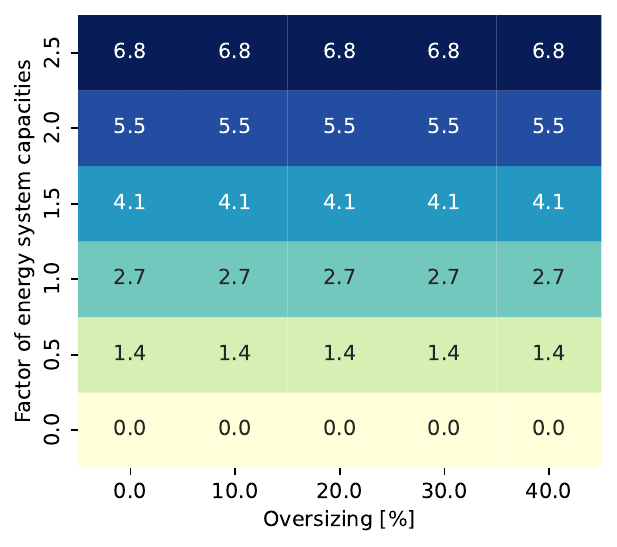}
         \caption{2022}
     \end{subfigure}
     \caption{TAC-optimal PV capacity [MWp].}
     \label{fig:supp:heatmap_PV}
\end{figure}

\begin{figure}[!ht]
     \begin{subfigure}[b]{0.32\textwidth}
         \centering
         \includegraphics[width=\textwidth]{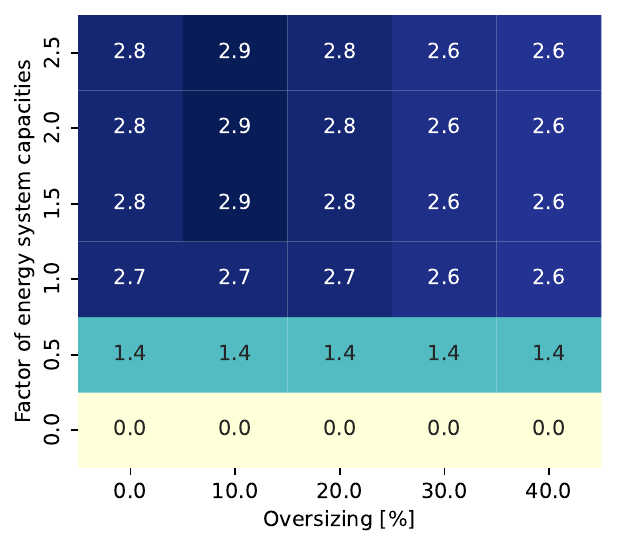}
         \caption{2020}
     \end{subfigure}
     \hfill
     \begin{subfigure}[b]{0.32\textwidth}
         \centering
         \includegraphics[width=\textwidth]{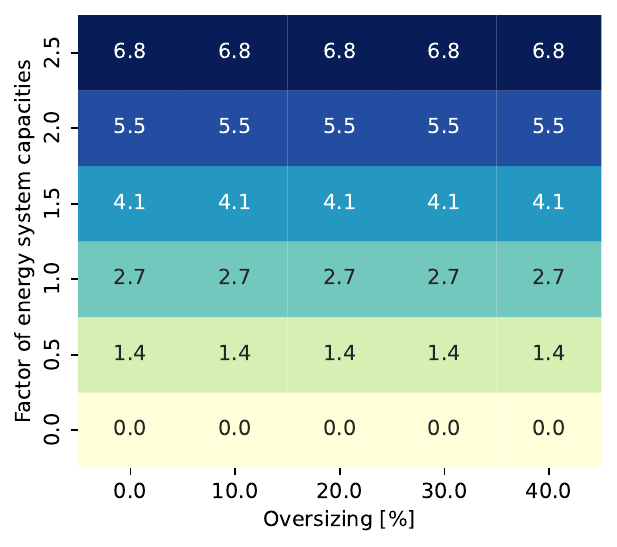}
         \caption{2021}
     \end{subfigure}
     \hfill
     \begin{subfigure}[b]{0.32\textwidth}
         \centering
         \includegraphics[width=\textwidth]{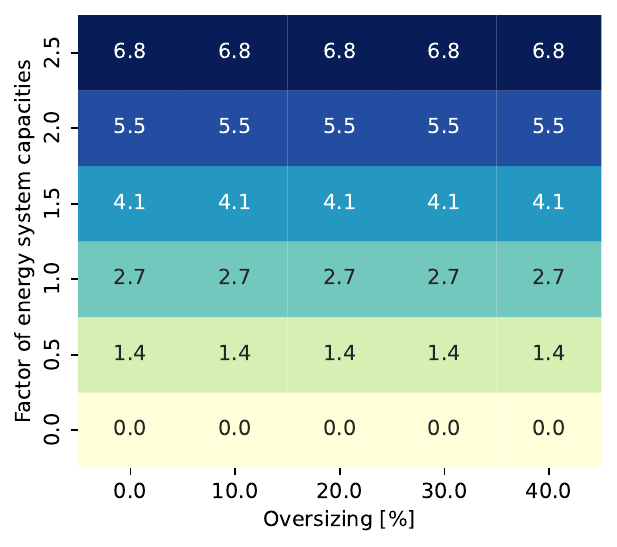}
         \caption{2022}
     \end{subfigure}
     \caption{TAC-optimal wind power capacity [MWp].}
     \label{fig:supp:heatmap_wind}
\end{figure}

\begin{figure}[!ht]
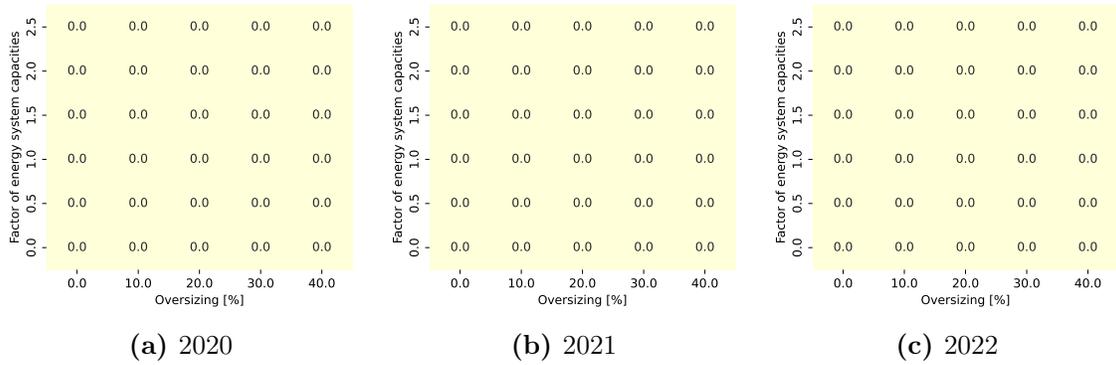

     \begin{subfigure}[b]{0.32\textwidth}
         \centering
         \includegraphics[width=\textwidth]{Figures/2020_Bat.pdf}
         \caption{2020}
     \end{subfigure}
     \hfill
     \begin{subfigure}[b]{0.32\textwidth}
         \centering
         \includegraphics[width=\textwidth]{Figures/2021_Bat.pdf}
         \caption{2021}
     \end{subfigure}
     \hfill
     \begin{subfigure}[b]{0.32\textwidth}
         \centering
         \includegraphics[width=\textwidth]{Figures/2022_Bat.pdf}
         \caption{2022}
     \end{subfigure}
     \caption{TAC-optimal battery capacity [MWh].}
     \label{fig:supp:heatmap_Bat}
\end{figure}

\begin{figure}
         \includegraphics[width=\textwidth]{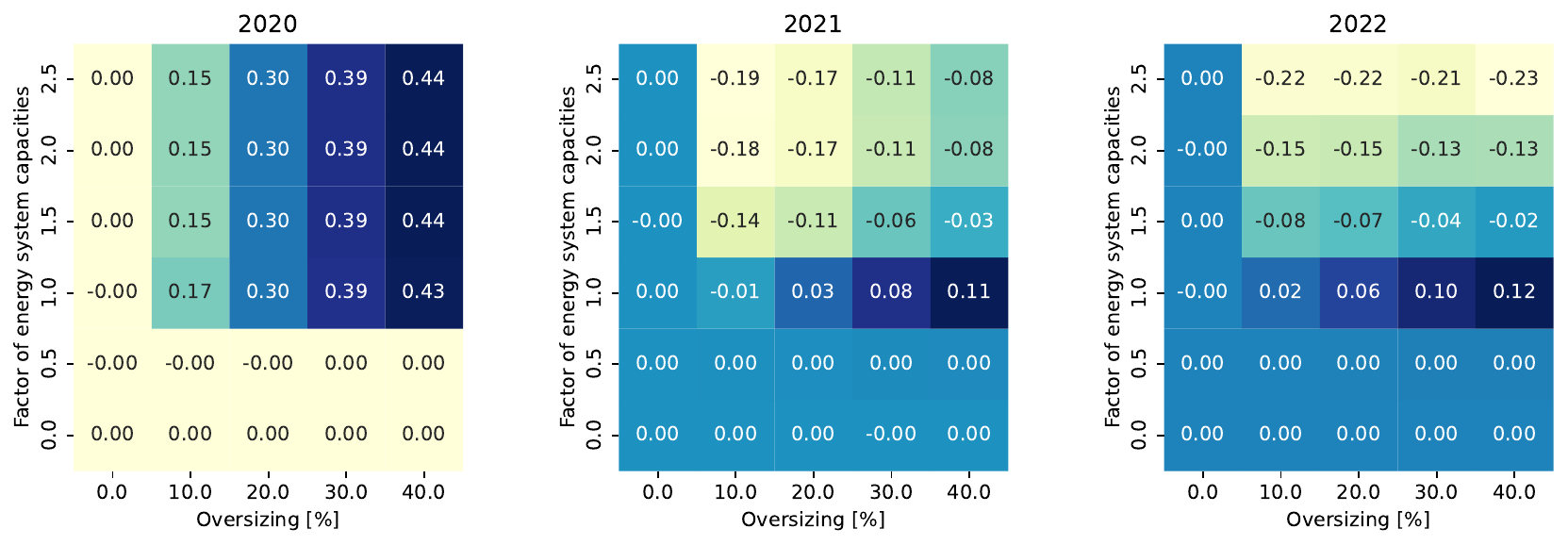}
         \caption{Relative difference between optimal TAC and estimated TAC [\%]: The estimated TAC is based on adding the savings from process flexibility and the installation of a local energy system. The difference is given relative to the optimal TAC shown in Fig. 7 in the main manuscript. }
         \label{fig:heatmap_diff_saves}
\end{figure}

Fig. 7 in the main manuscript shows the optimal TAC for a process with varying oversizing and varying maximum allowed energy system capacities. \cref{fig:supp:heatmap_PV,fig:supp:heatmap_Bat,fig:supp:heatmap_wind} show the corresponding optimal PV, wind power, and battery storage capacities. 
In 2020, PV is not built due to the low grid electricity price. In 2021, the higher electricity prices incentivize building PV up to 4.6 MWp. Finally, the even higher electricity prices in 2022 result in the maximum allowed PV capacities being installed. 
In 2020, wind power is incentivized up to 2.9 MWp. In contrast, the higher electricity prices in 2021 and 2022 result in maximum allowed wind power capacities being installed.
Battery storage capacity is not incentivized in any year.
Moreover, \cref{fig:heatmap_diff_saves} shows the relative difference between the optimal and estimated TAC.

\FloatBarrier

\section{Supporting material for the simultaneous market participation}

\begin{figure}[!ht]
        \centering
        \includegraphics[width=0.9\textwidth]{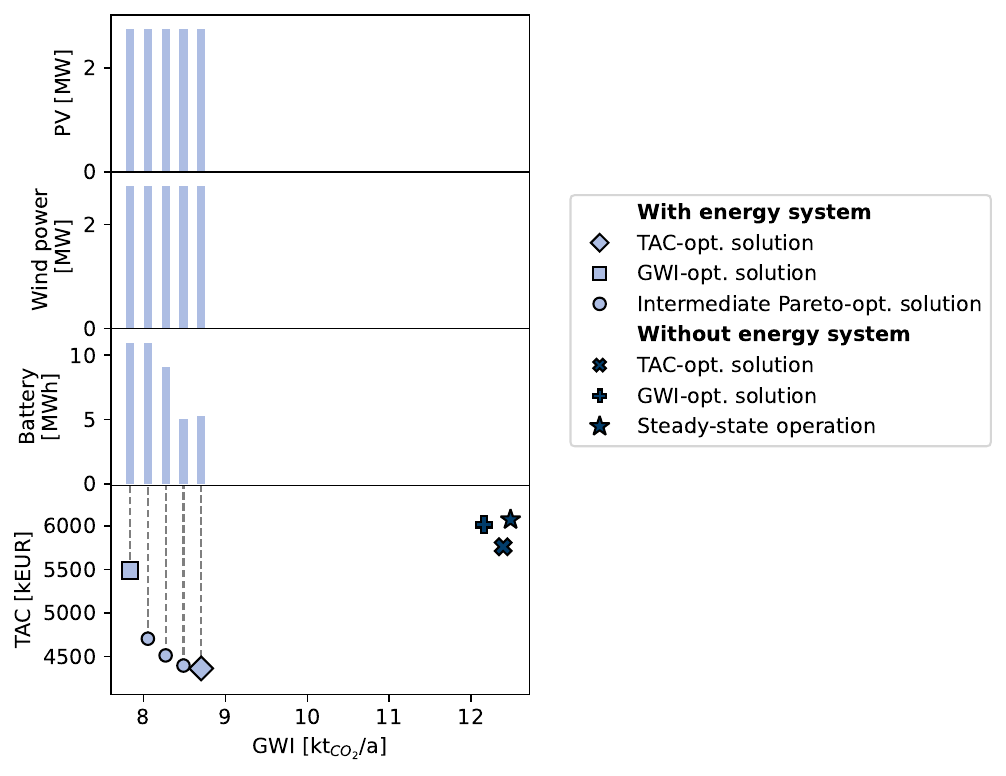}
         \caption{Energy system designs for 2022 for simultaneous \da{} and \id{} market participation: Five solutions are shown with the optimal capacities of the local energy system (upper three parts) and their GWI and TAC (lower part). TAC- and GWI-optimal \dr{} and steady-state operation without a local generation and storage system are given for comparison. }
         \label{fig:DAID_design}
\end{figure}

In Section 3.3 in the main manuscript, we discuss differences between simultaneous \da{} and \id{} market participation and \id{}-only participation. \cref{fig:DAID_design} shows the energy system designs for 2022 for the simultaneous market participation that differ from the \id{}-only participation (Fig. 4 (right) in the main manuscript) with respect to optimal battery capacities.

We assess the benefit of simultaneous market participation in Section 3.3 in the main manuscript based on 20 clusters. \cref{tab:DAID:cap10,tab:DAID:cap30,tab:DAID:saves10,tab:DAID:saves30} correspond to Tab. 6 and 7 in the main manuscript, but consider 10 and 30 clusters, respectively.

\cref{tab:DAID:cap10,tab:DAID:cap30} report the system capacities and show a minor variation of the optimal wind capacity in 2020 if 10 scenarios are used instead of 20 or 30. Recall that in Section 3.2 of the main manuscript, we show that the wind capacity can vary with respect to the process flexibility, i.e., it is rather sensitive with respect to the time series data of 2020.
Furthermore, \cref{tab:DAID:cap10,tab:DAID:cap30} show variations with respect to the trading amounts as these depend strongly on the electricity price realizations in the clusters.
For the year 2022, \cref{tab:DAID:cap10,tab:DAID:cap30} show that the optimal battery capacity increases with the number of clusters.
In \cref{tab:DAID:std}, we show that the average standard deviation for the price difference between \da {} and \id{} prices decreases with increasing number of clusters.
The decreased standard deviation allows for better positioning on the markets. 
Note that the battery capacity has a significant impact on the electricity trading amounts (\cref{tab:DAID:cap10,tab:DAID:cap30}), i.e., a larger battery enables larger trading amounts and vice versa.

\cref{tab:DAID:saves10,tab:DAID:saves30} show the TAC savings of the simultaneous market participation. Similar trends can be observed as those in Tab. 7 of the main manuscript. 
In comparison, both relative and absolute savings decrease and increase with less and more clusters, respectively. 
Again, we attribute this finding to the fact that more clusters lead to a lower standard deviation of the market deviation within a cluster and, thus, a better positioning on the markets.

\FloatBarrier

\begin{table}[!ht]
\caption{System capacities and electricity trading amounts based on 10 clusters}
    \label{tab:DAID:cap10}
    \centering
    \begin{tabular}{lrrr} \hline
 & 2020 & 2021 & 2022\\
 \hline  \multicolumn{3}{l}{\textbf{ID-only participation} } &  \\
\hline PV capacity & - & 2.74 MW & 2.74 MW \\
Wind capacity & 2.72 MW & 2.74 MW & 2.74 MW \\
Battery capacity & - & - & - \\
\hline Total purchases & 17,335 MWh & 15,958 MWh & 14,793 MWh \\
Total sales & 61 MWh & 216 MWh & 377 MWh \\
\hline \hline  \multicolumn{3}{l}{\textbf{Simultanous market participation} } &  \\
\hline PV capacity & - & 2.74 MW & 2.74 MW \\
Wind capacity & 2.73 MW & 2.74 MW & 2.74 MW \\
Battery capacity & - & - & 2.02 MWh \\
\hline DA purchases & 18,533 MWh & 16,777 MWh & 15,208 MWh \\
ID purchases & 6,808 MWh & 6,500 MWh & 9,965 MWh \\
Total purchases & 25,341 MWh & 23,277 MWh & 25,173 MWh \\
\hline DA sales & 156 MWh & 534 MWh & 2,859 MWh \\
ID sales & 7,936 MWh & 7,004 MWh & 7,749 MWh \\
Total sales & 8,092 MWh & 7,538 MWh & 10,608 MWh \\
\hline
\end{tabular}
\end{table}

\begin{table}[!ht]
\caption{System capacities and electricity trading amounts based on 30 clusters}
    \label{tab:DAID:cap30}
    \centering
    \begin{tabular}{lrrr} \hline
 & 2020 & 2021 & 2022\\
 \hline  \multicolumn{3}{l}{\textbf{ID-only participation} } &  \\
\hline PV capacity & - & 2.74 MW & 2.74 MW \\
Wind capacity & 2.74 MW & 2.74 MW & 2.74 MW \\
Battery capacity & - & - & - \\
\hline Total purchases & 17,311 MWh & 16,002 MWh & 14,802 MWh \\
Total sales & 60 MWh & 199 MWh & 364 MWh \\
\hline \hline  \multicolumn{3}{l}{\textbf{Simultanous market participation} } &  \\
\hline PV capacity & - & 2.74 MW & 2.74 MW \\
Wind capacity & 2.74 MW & 2.74 MW & 2.74 MW \\
Battery capacity & - & - & 6.3 MWh \\
\hline DA purchases & 17,020 MWh & 16,478 MWh & 19,307 MWh \\
ID purchases & 8,113 MWh & 7,564 MWh & 14,489 MWh \\
Total purchases & 25,134 MWh & 24,041 MWh & 33,796 MWh \\
\hline DA sales & 727 MWh & 1,328 MWh & 6,636 MWh \\
ID sales & 7,183 MWh & 6,918 MWh & 12,296 MWh \\
Total sales & 7,910 MWh & 8,246 MWh & 18,931 MWh \\
\hline
\end{tabular}
\end{table}

\begin{table}[!ht]
\caption{The TAC savings from simultaneous market participation based on 10 clusters}
    \label{tab:DAID:saves10}
    \centering
    \begin{tabular}{lrrr} \hline
 & 2020 & 2021 & 2022\\\hline
 \multicolumn{4}{l}{\textbf{Inflexible process without energy system} }  \\
Relative savings & 2.7\% & 1.5\% & 1.4\% \\
Absolute savings & 35.3 kEUR & 44.0 kEUR & 86.5 kEUR \\
\hline
 \multicolumn{4}{l}{\textbf{Flexible process without energy system} }  \\
Relative savings & 3.4\% & 1.9\% & 1.8\% \\
Absolute savings & 42.3 kEUR & 52.8 kEUR & 103.7 kEUR \\
\hline
  \multicolumn{4}{l}{\textbf{Flexible process with energy system} } \\
Relative savings & 3.4\% & 2.1\% & 2.5\% \\
Absolute savings & 42.6 kEUR & 54.0 kEUR & 115.6 kEUR \\
\hline
\end{tabular}    
\end{table}

\begin{table}[!ht]
\caption{The TAC savings from simultaneous market participation based on 30 clusters}
    \label{tab:DAID:saves30}
    \centering
    \begin{tabular}{lrrr} \hline
 & 2020 & 2021 & 2022\\\hline
 \multicolumn{4}{l}{\textbf{Inflexible process without energy system} }  \\
Relative savings & 3.5\% & 2.1\% & 2.2\% \\
Absolute savings & 46.1 kEUR & 59.8 kEUR & 135.7 kEUR \\
\hline
 \multicolumn{4}{l}{\textbf{Flexible process without energy system} }  \\
Relative savings & 4.4\% & 2.6\% & 2.8\% \\
Absolute savings & 55.4 kEUR & 71.8 kEUR & 162.8 kEUR \\
\hline
  \multicolumn{4}{l}{\textbf{Flexible process with energy system} } \\
Relative savings & 4.6\% & 3.1\% & 5.1\% \\
Absolute savings & 57.0 kEUR & 77.6 kEUR & 234.3 kEUR \\
\hline
\end{tabular}
\end{table}

\begin{table}[!t]
\caption{Average standard deviation of the market deviation for different numbers of clusters: The market deviation is the price difference between \da{} and \id{} prices. 
The clustering is determined as described in Section 2.4 of the main manuscript.
For each cluster, the standard deviation of each quarter-hourly time step is determined and the daily mean is taken. 
The listed values are the average of the daily mean values over all clusters.
The three years are considered separately.
}
    \label{tab:DAID:std}
    \centering
\begin{tabular}{lrrrr}
\hline
 &Unit & 10 clusters  & 20  clusters & 30 clusters  \\
\hline
2020       & [EUR/MWh]       & 13.0 & 11.6 & 10.7\\
2021        & [EUR/MWh]        & 22.9 & 20.6 & 19.9 \\
2022         & [EUR/MWh]       & 44.3 & 42.4 & 39.9 \\
\hline
\end{tabular}
\end{table}

\FloatBarrier

\bibliographystyle{apalike}

  \bibliography{BibTex.bib}